\documentclass[aos,MSNbibl,nameyear,dvips]{arximspdf}
\usepackage{dcolumn}
\usepackage{algpseudocode}
\usepackage{graphicx}

\doi{10.1214/14-AOS1245} %kopijuoti is PTS
\volume{42}
\issue{6}
\pubyear{2014}
\firstpage{2243}
\lastpage{2281}
\docsubty{FLA}

\makeatletter
\newcolumntype{d}[1]{D{.}{.}{#1}}
\newcommand{\rrvert}{\vert}
\newcommand{\llvert}{\vert}
\newcommand{\eqref}[1]{(\ref{#1})}

\long\def\comment#1{}

\newtheorem{theorem}{Theorem}[section]

\newtheorem{lemma}{Lemma}[section]

\newproclaim{defin}{Definition}[section]

\newproclaim{exercise}{Exercise}[section]

\newproclaim{assumption}{Assumption}[section]

\newproclaim{example}{Example}[section]

\newproclaim{remark}{Remark}[section]
\makeatother

\begin{document}
\begin{frontmatter}

\title{Wild binary segmentation for multiple change-point~detection}
\runtitle{Wild binary segmentation}

\begin{aug}
\author{\fnms{Piotr}~\snm{Fryzlewicz}\corref{}\ead[label=e1]{p.fryzlewicz@lse.ac.uk}}
\runauthor{P. Fryzlewicz}
\affiliation{London School of Economics}
\address{Department of Statistics\\
London School of Economics\\
Houghton Street\\
London WC2A 2AE\\
United Kingdom\\
\printead{e1}} %adresu isvedimo komanda gale!
\end{aug}

% HISTORY:
\received{\smonth{11} \syear{2013}}
\revised{\smonth{4} \syear{2014}}

% ABSTRACT
%
\begin{abstract}
We propose a new technique, called wild binary segmentation (WBS),
for consistent estimation of the number and locations of multiple
change-points in data. We assume that the number of change-points
can increase to infinity with the sample size.
Due to a certain random localisation mechanism, WBS works even for
very short spacings between the change-points and/or very small jump magnitudes,
unlike standard binary segmentation. On the other hand, despite its use of
localisation, WBS does not require the choice of a window or span parameter,
and does not lead to a significant increase in computational
complexity. WBS
is also easy to code. We propose two stopping criteria for WBS: one
based on thresholding and the other based on what we term the `strengthened
Schwarz information criterion'. We provide default recommended values
of the
parameters of the procedure and show that it offers very good practical
performance in comparison with the state of the art. The WBS methodology
is implemented in the R package \texttt{wbs}, available on CRAN.

In addition, we provide a new proof of consistency of binary segmentation
with improved rates of convergence, as well as a corresponding result for
WBS.
\end{abstract}

% KEYWORDS
\end{frontmatter}

%s1 #&#
\section{Introduction}\label{sec1}

A posteriori change-point detection problems have been of interest to
statisticians for many
decades. Although, naturally, details vary, a theme common to many of them
is as follows: a time-evolving quantity follows a certain stochastic
model whose
parameters are, exactly or approximately, piecewise constant. In such a
model, it is
of interest to detect the number of changes in the parameter values and
the locations
of the changes in time. Such piecewise-stationary modelling can be
appealing for a
number of reasons: the resulting model is usually much more flexible
than the corresponding
stationary model but still parametric if the number of change-points is fixed;
the estimated change-points are often `interpretable' in the sense that their
locations can be linked to the behaviour of some exogenous quantities
of interest; the
last estimated segment can be viewed as the `current' regime of
stationarity, which
can be useful in, for example, forecasting future values of the
observed process.
Finally, a posteriori segmentation can be a useful exploratory step in
the construction
of more complex models in which the piecewise constant variables are
themselves treated
as random and evolving according to a certain, perhaps Markovian, mechanism.

Arguably the simplest, `canonical' model with change-points is that of
the form
%
%e1 #&#
\begin{equation}
\label{eq:model} X_t = f_t + \varepsilon_t,
\qquad t = 1, \ldots, T,
\end{equation}
where $f_t$ is a deterministic, one-dimensional, piecewise-constant
signal with
change-points whose number $N$ and locations $\eta_1, \ldots, \eta_N$
are unknown.
The sequence $\varepsilon_t$ is random and such that $\mathbb
{E}(\varepsilon_t)$ is
exactly or approximately zero. In the simplest case $\varepsilon_t$ are
modelled as
i.i.d., but can also follow more complex time series models. The task
is to estimate
$N$ and $\eta_1, \ldots, \eta_N$ under various assumptions on~$N$, the
magnitudes
of the jumps and the minimum permitted distance between the
change-point locations.
Being univariate, model (\ref{eq:model}) excludes, for example, many
interesting time series
segmentation problems in which the process at hand is typically
parameterised by
more than one parameter in each segment. However, it still provides a
useful training ground for change-point
detection techniques in the sense that if a given method fails to
perform in the simple
model (\ref{eq:model}), it should not typically be expected to perform
well in more
complex settings.

There is considerable literature on a posteriori multiple change-point
detection in different
variants of model (\ref{eq:model}). \citet{ya89} consider least-squares
estimation
of $f_t$ in the case of a fixed $N$ (either known or unknown), under
the assumption
of $\varepsilon_t$ being i.i.d. In the case of a known $N$, they show
the consistency
of the estimated change-point locations with the rate of $O_P(1)$. They
also propose
a penalised least-squares estimator of $N$ in the case when it is
unknown but bounded.
In the Gaussian case, the Schwarz criterion is used to estimate an
unknown but bounded
$N$ in \citet{y88}, and a more general criterion that is also linear in
the number of change-points
appears in \citet{l95}. For an unknown but bounded $N$, \citet{lm00}
consider penalised least-squares
estimation, with a penalty linear in the number of change-points, and
show its consistency
for the number and locations of change-points for dependent $\varepsilon
_t$'s, including
the cases of strong mixing and long-range dependence; see also \citet{l99a} for a discussion
and some extensions of this result and \citet{l05a} for some practical
proposals regarding
the adaptive choice of the penalty parameter. For a fixed~$N$, \citet{pc06} propose a
likelihood criterion with a
penalty depending not only on the number, but also on the locations of
change-points, favouring
more uniformly-spread estimated change-points. For an unknown $N$, \citet{l05}
propose least-squares estimation with a penalty originating from the
model selection
approach of \citet{bm01} and show the least-squares consistency of the
resulting estimator
of $f_t$ (not of the estimated change-points themselves).
\citet{bklmw09} use the least-squares criterion with a linear penalty on
the number of
change-points and, under the assumption of a
finite but arbitrary~$N$, show various theoretical results including
analogues of those of \citet{ya89}. More general forms of Schwarz-like penalties
are studied, for example, in \citet{w08} and \citeauthor{c11} (\citeyear{c11,c14}).

Often, a major drawback of change-point estimators formulated as multivariate
optimisation problems, such as those based on penalised least-squares
or log-likelihood
fits, is their computational complexity, which is typically of order $O(T^2)$
[see, e.g., \citet{al89} and \citet{jsbaaggstt05}],
a~prohibitively slow speed for large datasets. \citet{kfe12} propose an
algorithm,
called PELT, that reduces the complexity to $O(T)$ but under the
assumption of
change-points being separated by time intervals drawn independently from
a probability distribution, a set-up under which considerations of
statistical consistency
are impossible due to these spacings being too short. \citet{r10}
proposes an alternative
`pruned dynamic programming' algorithm with the aim of reducing the
computational
effort, which, however, remains of order $O(T^2)$ in the worst case. Both
algorithms are revisited in the simulations section of this paper. An
interesting
approach to change-point detection, in the context of piecewise-stationary
AR time series models rather than in model (\ref{eq:model}), appears in
\citet{dlry06}: the minimum description length is used as the criterion
for segmentation, and it is minimised using a genetic algorithm to
reduce computational complexity.

A different route to reducing the computational complexity of the multiple
change-point detection problem is taken by \citet{hll10} who consider
the least-squares
criterion with a total variation penalty, which enables them to use the LARS
algorithm of \citet{ehjt04} to compute the solution in $O(N T \log(T))$ time.
For a known $N$ (only), they prove consistency of the resulting
estimated change-point locations
with near-optimal rates. We note, however, that the total variation penalty
is not an optimal one for change-point detection; see \citet{cf11}, who reiterate
an argument made earlier in \citet{bd93}. The total variation penalty
is also considered in the context of peak/trough detection by \citet{dk01},
who propose the `taut string' approach for fast computation, and
in the context of multiple change-point detection by \citet{r09} [as part
of the fused lasso penalty, proposed by \citet{tsrz05} and equivalent to
taut string
in model (\ref{eq:model})] and \citet{rw14}, who also point out that the
main result
in \citet{r09} is erroneous. On the other hand, \citet{w95} uses the
traditional fast
discrete wavelet transform to detect change-points.

An informative review of some multiple change-point detection methods
(in the context
of DNA segmentation, but applicable more widely) appears in \citet{bm98}.
\citet{changepoint.info} is an online repository of publications and software
related to change-point detection.

Binary segmentation (BS) is a generic technique for multiple
change-point detection
in which, initially, the entire dataset is searched for one
change-point, typically
via a CUSUM-like procedure. If and once a change-point is detected, the
data are
then split into two (hence the name `binary') subsegments, defined by
the detected
change-point. A similar search is then performed on either subsegment,
possibly resulting
in further splits. The recursion on a given segment continues until a
certain criterion is satisfied on it. Unlike estimators resulting from
multi-dimensional
optimisation of a certain global criterion, such as the least-squares
estimators reviewed
above, BS is a `greedy' procedure in the sense that it is performed
sequentially,
with each stage depending on the previous ones, which are never
re-visited. On the other
hand, each stage is particularly simple and involves one-, rather than
multi-dimensional
optimisation. To the best of our knowledge, the first work to propose
BS in a stochastic process
setting was \citet{v81}, who showed consistency of BS for the number and
locations
of change-points for a fixed $N$, with rates of convergence of the
estimators of
locations, under certain technical conditions on the norm of the
cumulative sum of
the process $X_t$, which in that work was assumed to be multivariate.
Testing for
change-points at each stage of the BS procedure was performed via a
simple CUSUM
test; however, the stopping criterion was not easy to compute in
practice due
to randomness in the previously detected change-points.
\citet{venkatraman1993} outlines an interesting proof of the consistency
of BS
for $N$ and for the change-point locations, even for $N$ increasing
with $T$, albeit
with sub-optimal rates for the locations.

Interestingly, BS in a setting similar to \citet{v81} (for a fixed $N$
and with $\varepsilon_t$
following a linear process), reappears in \citet{b97}, but without
references to the earlier works cited above. \citet{ccs11} provide a
proof of
consistency of BS for the number of change-points in the case of fixed $N$
and i.i.d. normal $\varepsilon_t$; however, links between their result
and the
analogous consistency results obtained in the above papers are not established.

We also note that BS has an interpretation
in terms of `unbalanced Haar' wavelets; see \citet{f07a}. BS is used for
univariate
time series segmentation in \citet{fsr11} and \citet{cf12}, and for
multivariate, possibly
high-dimensional time series segmentation in \citet{cf12a}.

The benefits of BS include low computational complexity [typically of
order $O(T\log T)$],
conceptual simplicity, and the fact that it is usually easy to code,
even in more
complex models than (\ref{eq:model}). \citet{kfe12} describe it as
`arguably the most
widely used change-point search method'. On the other hand, the fact
that each stage
of BS involves search for a \textit{single} change-point means that BS may
be unsuitable
for some functions containing multiple change-points in certain
configurations. Indeed,
in one of our side results of the paper, we show that BS is only
consistent when the
minimum spacing between any two adjacent change-points is of order
greater than
$T^{3/4}$ (even in the `easiest' case of jump magnitudes being bounded
away from zero), so relatively large.

In this work, we attempt to capitalise on the popularity and other
benefits of
BS and propose a multiple change-point detection procedure, termed wild binary
segmentation (WBS), which inherits the main
strengths of BS but attempts to eliminate its weaknesses. The main idea
is simple.
In the first stage, rather than using a global CUSUM statistic that
uses the entire
data sample $(X_1, X_2, \ldots, X_T)$, we randomly draw (hence the term
`wild') a number of
subsamples, that is, vectors $(X_s, X_{s+1}, \ldots, X_e)$, where $s$
and $e$ are integers such that $1 \le s < e \le T$,
and compute the CUSUM statistic on each subsample. We then maximise
each CUSUM, choose the
\textit{largest} maximiser over the entire collection of CUSUMs, and take
it to be the
first change-point candidate to be tested against a certain threshold.
If it is
considered to be significant, the same procedure is then repeated
recursively to the left and
to the right of it. The hope
is that even a relatively small number of random draws will contain a
particularly
`favourable' draw in which, for example, the
randomly drawn interval $(s, e)$ contains only one change-point,
sufficiently separated from both
$s$ and $e$: a set-up in which our CUSUM estimator of the change-point
location works particularly
well as it coincides with the maximum likelihood estimator (in the case
of $\varepsilon_t$
being i.i.d. Gaussian). We provide a lower bound for the number of
draws that guarantees
such favourable draws with a high probability. Apart from the threshold-based
stopping criterion for WBS, we also introduce another, based on what we
call the
\textit{strengthened Schwarz information criterion}.

By `localising' our CUSUM statistic in this randomised manner, we overcome
the issue of the `global' CUSUM being unsuitable for certain configurations
of multiple change-points. We also dramatically reduce the
permitted spacing between neighbouring change-points in comparison to
standard BS,
as well as the permitted jump magnitudes.
Moreover, by drawing intervals of different lengths, we avoid the
problem of
span or window selection, present in some existing approaches to
localising the CUSUM statistic,
for example in the `moving sum' (MOSUM) technique of \citet{hs01} and
\citet{km11}, and the (windowed) `circular' binary segmentation of
\citet{ovlw04}. We note that \citet{mj12} provide theoretical consistency
results for a method related to the latter, but not windowed and hence
computationally intensive, in the case of a bounded number of change-points.

The WBS procedure is computationally fast, consistent, as well as being
provably better than BS
and near-optimal in terms of the rates of convergence of the estimated
locations of
change-points even for very short spacings between neighbouring change-points
and for $N$ increasing with $T$. It also performs very well in practice and
is easy to code. Its R implementation is provided in the R package
\texttt{wbs}
[\citet{wbs}],
available from CRAN.

The paper is organised as follows. In Section~\ref{sec:mot}, we
motivate the
WBS procedure. In Section~\ref{sec:met}, we recall standard binary
segmentation (with some new consistency results) and outline the WBS technique
in more detail, also with corresponding results. In Section~\ref{sec:sim}, we give recommendations
on default parameter values and illustrate the performance of WBS in a
comparative
simulation study. In Section~\ref{sec:real}, we exhibit its performance
in the problem of
segmenting a time series arising in finance.

%s2 #&#
\section{Motivation}
\label{sec:mot}

In this work, we consider the model
%
%e2 #&#
\begin{equation}
\label{eq:mod} X_t = f_t + \varepsilon_t,
\qquad t = 1, \ldots, T,
\end{equation}
where $f_t$ is a deterministic, one-dimensional, piecewise-constant
signal with
change-points whose number $N$ and locations $\eta_1, \ldots, \eta_N$
are unknown.
Further technical assumptions on $f_t$ and $\varepsilon_t$ will be specified
later.

The basic ingredient of both the standard BS algorithm and WBS is the
CUSUM statistic defined by the inner product between the vector
$(X_s, \ldots,\break  X_e)$ and a particular vector of `contrast' weights
given below:
%
%e3 #&#
\begin{equation}
\label{eq:ip} \tilde{X}_{s,e}^b = \sqrt{\frac{e-b}{n(b-s+1)}}
\sum_{t=s}^b X_t - \sqrt {
\frac{b-s+1}{n(e-b)}}\sum_{t=b+1}^e
X_t,
\end{equation}
where $s \le b < e$, with $n = e-s+1$. It is used in different ways in
both algorithms. In its first step,\vspace*{-1pt}
the BS algorithm computes $\tilde{X}_{1,T}^b$ and then takes $b_{1,1} =
\arg\max_{b:1\le b<T} |\tilde{X}_{1,T}^b|$ to be the first
change-point candidate, whose significance is to be judged against a certain
criterion. If it is considered significant, the domain $[1, T]$ is
split into
two sub-intervals to the left and to the right of $b_{1,1}$ (hence the name
`binary segmentation'), and the recursion continues by computing
$\tilde{X}_{1,b_{1,1}}^b$ and $\tilde{X}_{b_{1,1}+1, T}^b$, possibly resulting
in further splits. The complete BS algorithm is outlined in Section~\ref{sec:sbs}.

We note that the maximisation of $|\tilde{X}_{s,e}^b|$ is equivalent to
the least
squares fit of a piecewise-constant function with one change-point to
$X_s^e = (X_s, \ldots,\break  X_e)'$,
in the following sense. Define ${\mathcal F}_{s,e}^b$ to be the set of
vectors supported
on $[s,e]$ with a single change-point at $b$. We have
\[
\arg\max_{b:s\le b<e}  \bigl|\tilde{X}_{s,e}^b\bigr | =
\arg\min_{b:s\le b<e} \min_{\bar{f}_{s,e}^b \in{\mathcal F}_{s,e}^b}
\bigl\|X_s^e - \bar{f}_{s,e}^b
\bigr\|_2^2.
\]
Therefore, if the true function $f_t$ contains only one change-point
$b_0$ on $[s,e]$,
then $\hat{b}_0 = \arg\max_{b:s\le b<e} |\tilde{X}_{s,e}^b|$ is the
least-squares estimator of $b_0$,
coinciding with the MLE in the case of $\varepsilon_t$ being i.i.d. Gaussian.
Speaking heuristically, this means that if $f_t$ contains only one
change-point on its
entire domain $[1, T]$, then $b_{1,1}$, the estimator of its location
from the first step
of the BS algorithm, is likely to perform well.

However, in the case of more than one change-point, the first step of
the BS algorithm
amounts to fitting $\bar{f}_{1,T}^b$, a function with a single
change-point, to data
with underlying multiple change-points, that is, to fitting the wrong
model. This may have
disastrous consequences, as the following example demonstrates.

%f1 #&#
\begin{figure}

\includegraphics{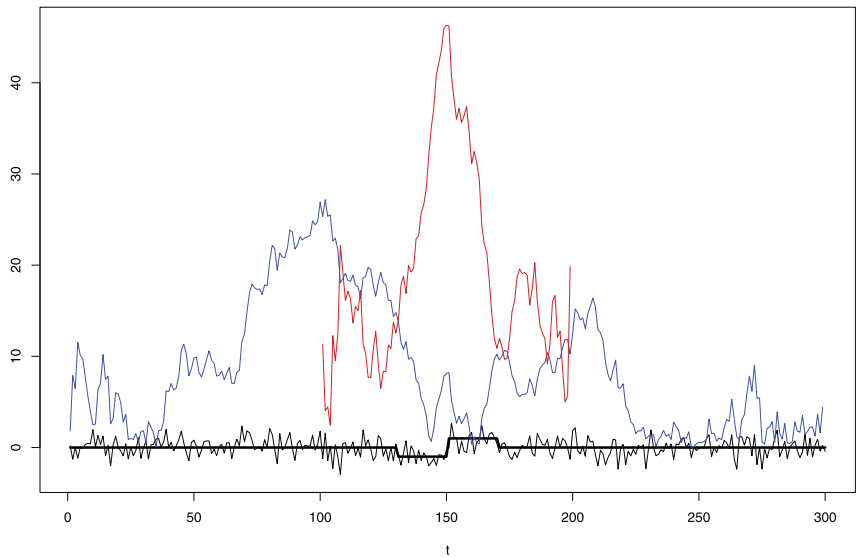}

\caption{True function $f_t$, $t = 1, \ldots, T=300$ (thick black),
observed $X_t$ (thin black),
$|\tilde{X}_{1,300}^b|$ plotted for $b = 1, \ldots,299$ (blue), and
$|\tilde{X}_{101,200}^b|$ plotted for $b = 101, \ldots, 199$
(red).}\label{fig:globloc}
\end{figure}

The function $\{f_t\}_{t=1}^{300}$ in Figure~\ref{fig:globloc} has
three change-points
(at $t=130, 150,\break  170$) which are concentrated in the middle of $f_t$,
and which `work
against each other' in the sense that the jump at $t=150$ is offset by
the two jumps
at $t=130, 170$. In the first step of BS, $|\tilde{X}_{1,300}^b|$ is
computed. However,
because of this unfavourable configuration of the change-points, its
maximum, occuring
around $b = 100$, completely misses all of them.

%f2 #&#
\begin{figure}

\includegraphics{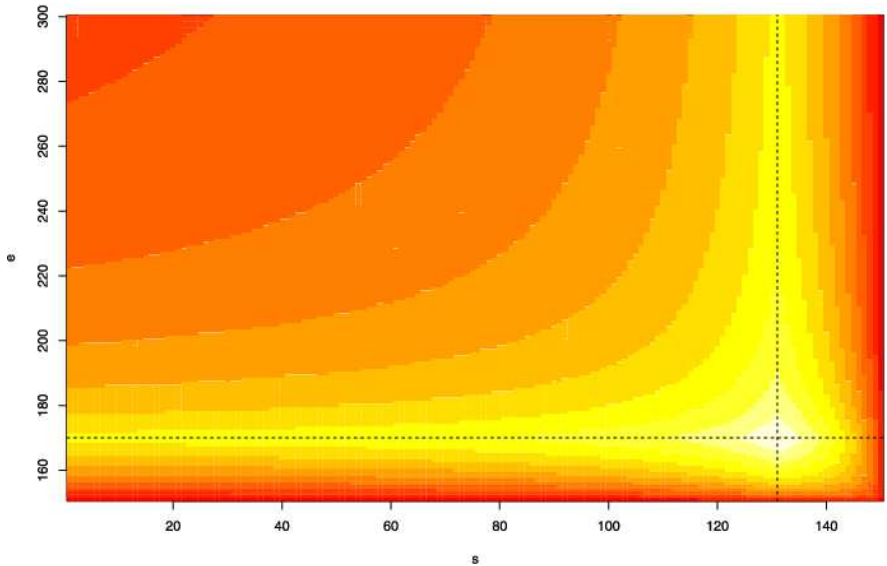}

\caption{Heat map of the values of $|\tilde{f}_{s,e}^{150}|$ as a
function of $s$ and $e$:
the lighter the colour, the higher the value. The two dashed lines
indicate the location
of the maximum, $(s,e)=(131,170)$.} \label{fig:ex}
\end{figure}

On the other hand, $|\tilde{X}_{101,200}^b|$ is successful in locating
the middle
change-point. Heuristically speaking, this is because the localised feature
(defined by the three change-points) is more `obvious' when
considered as part of the interval $[101, 200]$ than $[1, 300]$, in the sense
of the absolute inner product $|\tilde{f}_{101,200}^{150}|$ being much higher
than $|\tilde{f}_{1,300}^{150}|$ [where $\tilde{f}_{s,e}^b$ is defined
as in
(\ref{eq:ip}) but with $X$ replaced by~$f$]. This effect would be even more
pronounced if we `moved' the starting point of the inner product from
$s = 101$
towards the first change-point $t = 130$, and analogously the end point
$e = 200$
towards $t = 170$. In this example, the
inner product $|\tilde{f}_{s,e}^{150}|$ is maximised exactly when $s =
131$, $e = 170$
(i.e., when $s,e$ coincide with the two outside change-points), as this
creates the
`maximal' interval $[s,e]$ containing only the one change-point at $t = 150$.
This is further illustrated in Figure~\ref{fig:ex}.

Obviously, in practice, we cannot use the knowledge of the change-point
locations
to choose favourable locations for the start-point $s$ and the endpoint
$e$ of the inner product
$|\tilde{X}_{s,e}^b|$. We also cannot test all possible locations $s,
e$ as
this would be computationally prohibitive. Our main proposal
in this work is to randomly draw a number of pairs $(s,e)$ and
find $\arg\max_{b:s\le b<e}|\tilde{X}_{s,e}^b|$ for each draw. If the
number of draws is
suitably large, we will be able to guarantee, with high probability, a~particularly
favourable draw for which $[s,e]$ is long enough and only
contains one change-point at a sufficient distance from its endpoints
(or is sufficiently `close' to that situation, as in the example above).
The hope is that $\arg\max_{b:s\le b<e}|\tilde{X}_{s,e}^b|$
corresponding to that particular
draw will be a clear indicator of a true change-point in $f_t$. One perhaps
surprising aspect of this procedure is that the number of draws
guaranteed to achieve
this (for all change-points at once) is not large, as will be shown later.

This motivating discussion leads us to propose, in the next section,
the wild binary segmentation
algorithm for multiple change-point detection.

%s3 #&#
\section{Methodology and theory of wild binary segmentation}
\label{sec:met}

%s3.1 #&#
\subsection{Model and technical assumptions}

We make the following assumption.

%as3.1 #&#
\begin{assumption}
\label{ass:mainass}
(i) The random sequence $\{\varepsilon_t\}_{t=1}^T$ is i.i.d.
Gaussian with mean
zero and variance 1.

(ii) The sequence $\{f_t\}_{t=1}^T$ is bounded, that is, $|f_t|
< \bar{f} < \infty$ for $t = 1, \ldots, T$.
\end{assumption}

Assumption~\ref{ass:mainass}(i) is made both for technical convenience
and for clarity of exposition; it is reasonable to expect that it could in
principle be extended to dependent, heterogeneous and/or non-Gaussian noise.
We assume that $\operatorname{Var}(\varepsilon_t)$ is known, the reason being
that in practice it can
usually be estimated accurately using, for example, median absolute
deviation [\citet{h74}].
Such an assumption is standard in the literature on function estimation
in Gaussian
noise.

Different assumptions on the spacing between change-points and on the jump
magnitudes will be needed
by standard binary segmentation and by WBS. In what follows, denote
$\eta_0 = 0,
\eta_{N+1} = T$.

%as3.2 #&#
\begin{assumption}[(for standard binary segmentation)]
\label{ass:sbs}
The minimum spacing between change-points satisfies
$\min_{i=1, \ldots, N+1} |\eta_{i} - \eta_{i-1}| \ge\delta_T$, where
$\delta_T \ge C T^\Theta$ for $C > 0$, with $\Theta\le1$.
In addition, the magnitudes $f'_i = |f_{\eta_i} - f_{\eta_i - 1}|$
of the jumps
satisfy $\min_{i=1,\ldots,N} f'_i \ge\underline{f}_T$, where
$\underline{f}_T \ge C T^{-\varpi}$, with $\varpi\ge0$. The parameters
$\Theta$ and $\varpi$ satisfy $\Theta- \frac{\varpi}{2} > \frac{3}{4}$.
\end{assumption}

%as3.3 #&#
\begin{assumption}[(for WBS)]
\label{ass:wbs}
The minimum spacing between change-points satisfies
$\min_{i=1, \ldots, N+1} |\eta_{i} - \eta_{i-1}| \ge\delta_T$, and
the magnitudes $f'_i = |f_{\eta_i} - f_{\eta_i - 1}|$
of the jumps
satisfy $\min_{i=1,\ldots,N} f'_i \ge\underline{f}_T$, where
$\delta_T$ and $\underline{f}_T$ are linked by the requirement
$\delta_T^{1/2} \underline{f}_T \ge C \log^{1/2}T$ for a large
enough $C$.
\end{assumption}

It is worth noting that we do not assume any further upper bounds on the
number $N$ of change-points, other than those implied by the minimum spacing
$\delta_T$. In other words, $N$ can be as large as allowed by $\delta_T$,
and in particular can increase to infinity with $T$. Therefore, formally,
we have $N = N(T)$ and $\eta_i = \eta_i(T)$ for $i = 1, \ldots, N+1$. However,
for economy of notation and keeping in line with many other papers on
change-point
detection, in the remainder of the paper we use the shorthand notation
$N, \eta_i$ rather than the longer notation $N(T), \eta_i(T)$.

The quantity $\delta_T^{1/2} \underline{f}_T$ appearing in Assumption~\ref{ass:wbs} is well known in the `statistical signal detection'
literature. For example, \citet{cw13} summarise results which show that
detection of hat-shaped signals observed in Gaussian noise is impossible
if (the equivalent of) this quantity is below a certain threshold. See
also \citet{ds01} and \citet{fms14} for related discussions. We will argue
in Section~\ref{th:wbs} that our Assumption~\ref{ass:wbs} is
rate-near-optimal from
this point of view.

%s3.2 #&#
\subsection{Standard binary segmentation}
\label{sec:sbs}

To gain a better understanding of the improvement offered by WBS over
standard BS, we first provide a theoretical consistency
result for the latter. The BS algorithm is best defined recursively and
hence described by pseudocode. The main function is defined as follows.

\begin{algorithmic}
\Function{BinSeg}{$s$, $e$, $\zeta_T$}
\If{$e-s < 1$}
\State STOP
\Else
\State$b_0:= \arg\max_{b\in\{s, \ldots, e-1\}} |\tilde{X}_{s,e}^b|$
\If{$|\tilde{X}_{s,e}^{b_0}| > \zeta_T$}
\State add $b_0$ to the set of estimated change-points
\State\textsc{BinSeg}($s$, $b_0$, $\zeta_T$)
\State\textsc{BinSeg}($b_0+1$, $e$, $\zeta_T$)
\Else
\State STOP
\EndIf
\EndIf
\EndFunction
\end{algorithmic}

Given the above definition, the standard BS procedure is launched
by the call \textsc{BinSeg}($1$, $T$, $\zeta_T$), where $\zeta_T$ is a threshold
parameter. Let $\hat{N}$ denote the number of change-points estimated by
the BS algorithm, and $\hat{\eta}_1, \ldots, \hat{\eta}_{\hat{N}}$ their
locations, sorted in increasing order. The following consistency
theorem holds.

%th3.1 #&#
\begin{theorem}
\label{th:sbs}
Let $X_t$ follow model (\ref{eq:mod}), and
suppose Assumptions \ref{ass:mainass} and~\ref{ass:sbs} hold. Let $N$ and
$\eta_1, \ldots, \eta_N$ denote, respectively, the number and locations of
change-points. Let $\hat{N}$ denote the number, and $\hat{\eta}_1,
\ldots, \hat{\eta}_{\hat{N}}$
the locations, sorted in increasing order, of the change-point estimates
obtained by the standard binary segmentation algorithm. Let the threshold
parameter satisfy $\zeta_T = c_1 T^\theta$ where $\theta\in(1-\Theta, \Theta-1/2-\varpi)$
if $\Theta\in(\frac{3}{4}, 1)$, or $\zeta_T \ge c_2 \log^p T$ ($p >
1/2$) and $\zeta_T \le c_3 T^\theta$
($\theta< 1/2 - \varpi$) if $\Theta= 1$, for any positive constants
$c_1, c_2, c_3$. Then there
exist positive constants $C$, $C_1$ such that $P(\mathcal{A}_T) \ge1 -
C_1T^{-1}$, where
\[
\mathcal{A}_T = \Bigl\{ \hat{N} = N; \max_{i=1,\ldots,N}|
\hat{\eta}_i - \eta_i| \le C\epsilon_T \Bigr
\}
\]
with $\epsilon_T = T^2 \delta_T^{-2} (\underline{f}_T)^{-2} \log T$.
\end{theorem}

We note that the rates of convergence of $\hat{\eta}_i$ are better than
those obtained by \citet{venkatraman1993} and \citet{fsr11},
both of which consider consistency of the BS procedure for the
number of change-points $N$ possibly increasing with $T$; they are also
better than those in \citet{cf12} (where $N$ is assumed to be bounded).
The latter three papers use the assumption that $\underline{f}_T$ is
bounded away from zero.
The improvement is due to the crucial and new Lemma~\ref{lem:vnear}. Rates are particularly important here, as they inform the
stopping criterion (i.e., the admissible magnitude of the threshold
$\zeta_T$),
rather than merely quantifying the performance of the procedure.

As an aside, we mention that in the case $\delta_T = o(T)$, it is
possible to further improve our
rates via a simple trick, whereby change-point locations are re-estimated
by maximising the CUSUM statistic $|\tilde{X}_{s,e}^b|$ on each interval
$[s,e]$ where $s,e$ are respective mid-points of two adjacent intervals
$[\hat{\eta}_{i-1}+1, \hat{\eta}_{i}]$, $[\hat{\eta}_{i}+1, \hat{\eta}_{i+1}]$
(with the convention $\hat{\eta}_0 = 0, \hat{\eta}_{\hat{N}+1} = T$).
This refinement can be applied to any multiple
change-point detection procedure, not just BS. However, even with this
refinement, the BS procedure as defined above is only guaranteed to produce
valid results under Assumption~\ref{ass:sbs}, which is rather restrictive
in terms of the permitted distance between change-points and the magnitudes
of the jumps.

%s3.3 #&#
\subsection{Wild binary segmentation}
\label{sec:wbssub}

Denote by $F_T^M$ a set of $M$ random intervals
$[s_m, e_m]$, $m = 1, \ldots, M$, whose start- and end-points
have been drawn (independently with replacement) uniformly
from the set $\{1, \ldots, T\}$. Guidance as to a suitable
choice of $M$ will be given later. Again using pseudocode, the main function
of the WBS algorithm is defined as follows.

\begin{algorithmic}
\Function{WildBinSeg}{$s$, $e$, $\zeta_T$}
\If{$e-s < 1$}
\State STOP
\Else
\State${\mathcal M}_{s,e}:= $ set of those indices $m$ for which
$[s_m, e_m] \in F_T^M$ is such that
\hspace*{32pt}$[s_m,e_m]\subseteq[s,e]$
\State(Optional: augment ${\mathcal M}_{s,e}:= {\mathcal M}_{s,e}
\cup\{0\}$, where $[s_0,e_0] = [s,e]$)
\State$(m_0, b_0):= \arg\max_{m\in{\mathcal M}_{s,e}, b\in\{s_m,
\ldots, e_m-1\}} |\tilde{X}_{s_m,e_m}^b|$
\If{$|\tilde{X}_{s_{m_0},e_{m_0}}^{b_0}| > \zeta_T$}
\State add $b_0$ to the set of estimated change-points
\State\textsc{WildBinSeg}($s$, $b_0$, $\zeta_T$)
\State\textsc{WildBinSeg}($b_0+1$, $e$, $\zeta_T$)
\Else
\State STOP
\EndIf
\EndIf
\EndFunction
\end{algorithmic}

The WBS procedure is launched by the call \textsc{WildBinSeg}($1$, $T$,
$\zeta_T$).
We believe that the WBS procedure is not difficult to code even for the
nonexpert, unlike some change-point detection algorithms based on
dynamic programming.
Let $\hat{N}$ denote the number of change-points estimated by
the WBS procedure, and $\hat{\eta}_1, \ldots, \hat{\eta}_{\hat{N}}$ their
locations, sorted in increasing order.

The optional augmentation of ${\mathcal M}_{s,e}$ by $\{0\}$ is done to ensure
that the algorithm also examines the entire current interval $[s,e]$,
and not
only its randomly drawn subintervals, in case $[s,e]$
only contains one change-point and hence it is optimal to examine
$[s,e]$ in its entirety.
We note that unlike the BS procedure, the WBS algorithm (in the
case without the optional augmentation) returns estimated change-points
in the order corresponding to decreasing maxima of $|\tilde{X}_{s_m,e_m}^b|$,
which is due to the maximisation over~$m$. There
is no corresponding maximisation in the BS procedure, which
means that the maxima of the CUSUM statistics corresponding
to estimated change-points in the latter procedure are not necessarily
arranged in decreasing order.

Finally, we motivate the use of random, rather than fixed, intervals.
As demonstrated in Section~\ref{sec:mot}, some change-points
require narrow intervals $[s,e]$ around them in order to be
detectable. For such change-points, the use of random intervals, as in
the WBS algorithm, means that there is always a positive probability,
sometimes high, of there being a suitably narrow interval around
them in the set $F_T^M$. On the other hand, consider a fixed design,
where the start-points $s_m$ and end-points $e_m$ take all possible
values from a fixed subset of $\{1, \ldots, T\}$, of such cardinality
that the number of resulting intervals is the same as in the random
design. For such a fixed design (however it is chosen), at least some
of the intervals will inevitably be significantly longer than the
corresponding random ones, so that they may not permit detection of
such change-points if those happen to
lie within them. Another reason is that through the use of randomness,
we avoid having to make the subjective choice of a particular fixed design.
Finally, if the number of intervals drawn turns out to be insufficient, it
is particularly easy to add further intervals if the design is random;
this is achieved simply by drawing further intervals from the same distribution.
In the case of a fixed design, the entire collection may need to be re-drawn
if the distribution of interval lengths is to be preserved.
However, for a very large number $M$ of intervals, the difference in performance
between the random and deterministic designs is likely to be minimal.

The following theoretical result holds for the WBS algorithm.

%th3.2 #&#
\begin{theorem}
\label{th:wbs}
Let $X_t$ follow model (\ref{eq:mod}), and
suppose Assumptions \ref{ass:mainass} and~\ref{ass:wbs} hold. Let $N$ and
$\eta_1, \ldots, \eta_N$ denote, respectively, the number and locations of
change-points. Let $\hat{N}$ denote the number, and $\hat{\eta}_1,
\ldots, \hat{\eta}_{\hat{N}}$
the locations, sorted in increasing order, of the change-point estimates
obtained by the wild binary segmentation algorithm.
There exist two constants $\underline{C}$, $\overline{C}$ such that if
$\underline{C} \log^{1/2}T \le\zeta_T \le\overline{C} \delta_T^{1/2}
\underline{f}_T$,
then $P(\mathcal{A}_T) \ge1 - C_1 T^{-1} - T \delta_T^{-1} (1 - \delta
_T^2 T^{-2} / 9)^M$, where
\[
\mathcal{A}_T = \Bigl\{ \hat{N} = N; \max_{i=1,\ldots,N}|
\hat{\eta}_i - \eta_i| \le C\log T (
\underline{f}_T)^{-2} \Bigr\}
\]
for certain positive $C$, $C_1$.
\end{theorem}

Some remarks are in order. Firstly, we note that Assumption~\ref
{ass:wbs} is
much milder than Assumption~\ref{ass:sbs}. As an illustration, consider
the case
when $\underline{f}_T$ is bounded away from zero (although we emphasise
that both
algorithms permit $\underline{f}_T \to0$, albeit at different rates).
In this case,
the WBS method produces consistent results even if the minimum spacing
$\delta_T$ between the true change-points is logarithmic in $T$,
whereas $\delta_T$
must be larger than $O(T^{3/4})$ in standard BS. Furthermore, for a
given separation
$\delta_T$ and minimum jump height $\underline{f}_T$, the admissible
range of
threshold rates for the WBS method is always larger than that for BS.
In this sense, the
WBS method may be viewed as more robust than BS to the possible
misspecification of the value of the threshold.

Secondly, unlike the BS algorithm, the lower bound for the threshold
$\zeta_T$
in the WBS method is always square-root logarithmic in $T$,
irrespective of the
spacing~$\delta_T$. This is also the only threshold rate that yields
consistency for any
admissible separation $\delta_T$ and minimum jump size $\underline
{f}_T$. For this reason, we use the rate $\log^{1/2} T$
as the default rate for the magnitude of the threshold, and hence, in
the remainder
of the article, we consider thresholds of the form $\zeta_T = C \sqrt
{2} \log^{1/2} T$
(we introduce the factor of $\sqrt{2}$ in order to facilitate the
comparison of $\zeta_T$
to the `universal' threshold in the wavelet thresholding literature,
which is of the form
$\sqrt{2} \log^{1/2} T$). Practical choice of the constant $C$ will be
discussed in Section~\ref{sec:sim}.
In BS, the only threshold rate that leads to consistency for any
admissible $\delta_T$
is $\zeta_T \sim T^{1/4 - \varpi/2}$ (where $\sim$ means `of the order of'
throughout the paper).

Thirdly, again unlike the BS algorithm, the rate of convergence of the
estimated change-point locations in the WBS method does not depend on the
spacing $\delta_T$ (as long as $\delta_T^{1/2}\underline{f}_T$ is large
enough in the sense of Assumption~\ref{ass:wbs}) but only on the minimum
jump height $\underline{f}_T$. We now consider the special case of
$\underline{f}_T$
being bounded away from zero, and discuss the optimality, up to at most
a logarithmic
factor, of wild binary segmentation in estimating the change-point
locations in this
setting. In the case $\delta_T \sim T$, the optimal rate in detecting
change-point
locations is $O_P(1)$ in the sense that for any estimator $\hat{\eta
}_i$ of $\eta_i$,
we have $|\hat{\eta}_i - \eta_i| = O_P(1)$ at best; see, for example,
\citet{k87}. This can be reformulated as
$P(|\hat{\eta}_i - \eta_i| \ge a_T) \to0$ for any sequence $a_T \to
\infty$.
In the case $\underline{f}_T > \underline{f} > 0$, the result of
Theorem~\ref{th:wbs}
implies $P(\exists_i  |\hat{\eta}_i - \eta_i| \ge C\log T) \to0$,
thus matching the above minimax result up to a logarithmic term.
However, we emphasise that
this is in the (more challenging) context where (i) the number $N$ of
change-points is possibly
unbounded with $T$, and (ii) the spacing $\delta_T$ between
change-points can be much shorter than
of order~$T$.

We now discuss the issue of the minimum number $M$ of random draws
needed to ensure that the bound on the speed of convergence of
$P(\mathcal{A}_T)$ to 1 in Theorem~\ref{th:wbs} is suitably small. Suppose
that we wish to ensure
\[
T \delta_T^{-1} \bigl(1 - \delta_T^2
T^{-2} / 9\bigr)^M \le T^{-1}
\]
in order to match the rate of the term $C_1 T^{-1}$ in the upper bound
for $1 - P(\mathcal{A}_T)$ in Theorem~\ref{th:wbs}.
Bearing in mind that $\log(1-y) \approx-y$ around
$y=0$, this is, after simple algebra, (practically) equivalent to
\[
M \ge\frac{9 T^2}{\delta_T^2} \log\bigl(T^{2} \delta_T^{-1}
\bigr).
\]
In the `easiest' case $\delta_T \sim T$, this results in a logarithmic
number of draws, which leads to particularly low computational complexity.
Naturally, the required $M$ progressively increases as $\delta_T$
decreases. Our practical recommendations for the choice of $M$ are discussed
in Section~\ref{sec:sim}.

Furthermore, we explain why the binary recursion is needed in the WBS algorithm
at all: the careful reader may wonder why change-points are not
estimated simply
by taking \textit{all} those points that attain the maxima of $|\tilde
{X}_{s_m, e_m}^b|$
exceeding the threshold $\zeta_T$, for all intervals $[s_m,e_m] \in
F_T^M$. This is
because such a procedure would very likely lead to some true
change-points being
estimated more than once at different locations. By proceeding sequentially
as in the WBS algorithm, and by restricting ourselves to those intervals
$[s_m, e_m]$ that fully fall within the current interval of interest $[s,e]$,
we ensure that this problem does not arise. Another reason for
proceeding sequentially is the optional augmentation of
${\mathcal M}_{s,e}$ by $\{0\}$ in the WBS algorithm, which depends on
the previously
detected change-points and hence is not feasible in a nonsequential
setting.

Regarding the optimality of the lowest permitted rate for $\delta
_T^{1/2}\underline{f}_T$
in Assumption~\ref{ass:wbs}, recall that, by Theorem~\ref{th:wbs},
$\delta_T$ must be at least as large as $\max_{i=1,\ldots,N}|\hat{\eta
}_i - \eta_i|$,
or it would not be possible to match the estimated change-point
locations with the true
ones. Therefore, $\delta_T$ cannot be of a smaller order than $\log T$.
By the minimax arguments summarised in \citet{cw13} (but using our notation),
the rate of the smallest possible $\delta_T^{1/2}\underline{f}_T$
that permits change-point detection (by any method) for this range of
$\delta_T$
is $(\log T - \log\log T)^{1/2}$. Our Assumption~\ref{ass:sbs}
achieves this rate
up to the negligible double-logarithmic factor and therefore is optimal
under the
circumstances.

Randomised methods are not commonly used in nonparametric statistics
(indeed, we are not aware of any other commonly used such method); however,
randomised techniques are beginning to make headway in statistics in
the context
of `big data'; see, for example, the review articles \citet{m10} and
\citet{hmt11}. The proof technique in Theorem~\ref{th:wbs}
relies on some subtle arguments regarding the guarantees of quality of
the randomly
drawn intervals.

%s3.4 #&#
\subsection{Strengthened Schwarz information criterion for WBS}
\label{sec:ssic}

Naturally, the estimated number $\hat{N}$ and locations $\hat{\eta}_1,
\ldots, \hat{\eta}_{\hat{N}}$
of change-points depend on the selected threshold $\zeta_T$. For the
purpose of this paragraph,
denote $\hat{N}(\zeta_T) = \hat{N}$ and $\mathcal{C}(\zeta_T) = \{\hat
{\eta}_1, \ldots, \hat{\eta}_{\hat{N}(\zeta_T)}\}$.
It is a property of the WBS method that $\hat{N}(\zeta_T)$ is a nondecreasing
function of $\zeta_T$, each increase has size 1 almost-surely, and the
collection $\mathcal{C}(\zeta_T)$
is nested in the sense that if $\zeta'_T < \zeta_T^{\prime\prime}$ then
$\mathcal{C}(\zeta_T^{\prime\prime}) \subseteq\mathcal{C}(\zeta'_T)$. Consider
any decreasing sequence $\{\zeta_T^k\}_{k=0}^K$ of thresholds such that
$|\mathcal{C}(\zeta_T^k)| = k$ for a certain fixed constant $K$, and assume
$N \le K$. One may perform model selection either by choosing a
suitable threshold $\zeta_T$ and hence selecting the associated model
$\mathcal{C}(\zeta_T)$, or alternatively by considering the sequence
of model candidates $\{\mathcal{C}(\zeta_T^k)\}_{k=0}^K$ and choosing
one that
optimises a certain criterion, thereby by-passing the question of threshold
choice entirely. Thus it is a viable alternative to view the
`solution path' $\mathcal{C}(\zeta_T^k)$ not as a
function of threshold $\zeta_T^k$, but as a function of the number $k$
of change-point candidates. We define
$\mathcal{C}_{k} = \mathcal{C}(\zeta_T^k)$.
In this section, we propose to select a model out of the collection
$\{\mathcal{C}_{k}\}_{k=0}^K$ by minimising what
we term the `strengthened Schwarz information criterion' (sSIC),
defined as follows.

For any candidate model $\mathcal{C}_{k}$, denote by $\hat{f}^k_t$ the estimate
of $f_t$ defined by $\hat{f}^k_t = (\hat{\eta}_{i+1}-\hat{\eta}_i)^{-1}
\sum_{j=\hat{\eta}_i+1}^{\hat{\eta}_{i+1}} X_j$
for $\hat{\eta}_i+1 \le t \le\hat{\eta}_{i+1}$. Let $\hat{\sigma}^2_k
= T^{-1} \sum_{t=1}^T (X_t - \hat{f}^k_t)^2$
be the corresponding maximum likelihood estimator of the residual
variance. We define
%
%e4 #&#
\begin{equation}
\label{eq:ssic} \operatorname{sSIC}(k) = \frac{T}{2} \log \hat{
\sigma}^2_k + k \log^\alpha T.
\end{equation}
We remark that the choice $\alpha= 1$ corresponds to the standard SIC
penalty, considered, for example, by \citet{y88} in the context of
multiple change-point detection in a model similar to
ours performed via a full penalised least-squares minimisation. The
following result
holds.

%th3.3 #&#
\begin{theorem}
\label{th:ssic}
Let $X_t$ follow model (\ref{eq:mod}), and let the assumptions of
Theorem~\ref{th:wbs} hold.
Let $N$ and $\eta_1, \ldots, \eta_N$ denote, respectively, the number
and locations of
change-points. Let $N \le K$, where $K$ is a certain constant
independent of $T$.
Let the constant $\alpha> 1$ be such that $\log^\alpha T = o(\delta_T
\underline{f}_T^2)$. Let the candidate
models $\{\mathcal{C}_k\}_{k=0}^K$ be produced by the WBS algorithm,
and let
$\hat{N} = \arg\min_{k=0, \ldots, K} \operatorname{sSIC}(k)$.
Then $P(\mathcal{A}_T) \ge1 - C_1 T^{-1} - T \delta_T^{-1} (1 - \delta
_T^2 T^{-2} / 9)^M$, where
\[
\mathcal{A}_T = \Bigl\{ \hat{N} = N; \max_{i=1,\ldots,N}|
\hat{\eta}_i - \eta_i| \le C\log T (
\underline{f}_T)^{-2}\Bigr\}
\]
for certain positive $C$, $C_1$.
\end{theorem}

The only parameter of the above procedure is the constant
$\alpha$, and we require that $\alpha> 1$, which results in a stronger
penalty than in the standard SIC, hence the term `strengthened' SIC.
Noting the requirement that $\log^\alpha T = o(\delta_T \underline
{f}_T^2)$, we focus attention on
values of $\alpha$ close to 1, to ensure the admissibility of the
sSIC criterion for as large a class of signals as possible; from this
point of
view, it is tempting to regard this region of the parameter space for
$\alpha$
as a natural default choice. With this in mind, in the remainder of the paper,
we report the performance of sSIC with $\alpha= 1.01$, which also
ensures that
the results remain close to those obtained by SIC.

We further note that unlike in thresholding, where the magnitude
of the threshold is sensitive to $\operatorname{Var}(\varepsilon_t)$, the
minimisation of the
sSIC penalty in (\ref{eq:ssic}) is independent of
$\operatorname{Var}(\varepsilon_t)$ due to the use of the logarithmic
transformation in
$\log  \hat{\sigma}^2_k$. This logarithmic transformation causes
$\operatorname{Var}(\varepsilon_t)$ to have an additive contribution to the sSIC
criterion in (\ref{eq:ssic}), and therefore this term has no impact on the
minimisation.

In summary, the attraction of the sSIC approach lies in the fact that
the default choice of the parameter of the procedure is perhaps easier than
in the thresholding approach. On the other hand, the theoretical validity
of sSIC in the version of Theorem~\ref{th:ssic} requires that $N \le K$ for
a finite $K$ and that the lowest admissible $\delta_T \underline
{f}_T^2$ is (marginally) larger than in the thresholding
approach. The requirement of a finite $K$
is common to penalised approaches to multiple change-point detection;
see, for example,
\citet{y88} and \citet{c14}.\vadjust{\goodbreak}

%s4 #&#
\section{Parameter choice and simulation study}
\label{sec:sim}

%s4.1 #&#
\subsection{Parameter choice}

We now elaborate on the choice of the number $M$ of the random
draws, and the threshold constant $C$.

\textit{Choice of $M$.}
The parameter $M$ should be chosen to be `as large as possible'
subject to computational constraints. We note that with the
optional augmentation of ${\mathcal M}_{s,e}$ by $\{0\}$,
the WBS reduces to standard BS for $M = 0$, so even a relatively
small value of $M$ is likely to bring benefits in terms of
performance. Our recommendation is to set $M = 5000$ for datasets
of length $T$ not exceeding a few thousand. As an
example, with this value of $M$, we achieved the average
computation time of 1.20 seconds for a dataset of length
$T = 2000$. The code was written in a combination of R and C, and
executed on a 3.40~GHz quad-core with 8~GB of RAM, running Windows 7.
The implementation of WBS in the R package \texttt{wbs} is faster still.

Moreover, the larger the value of $M$, the more negligible the
dependence of the
solution on the particular random draw. For $M = 5000$, this dependence has
been observed to be very minimal.

\textit{Choice of the threshold constant $C$.}
In Section~\ref{sec:wbssub}, we motivate the use of thresholds of
the form $\zeta_T = C \sqrt{2} \log^{1/2} T$. There remains the
question of how to
choose the threshold constant $C$. We firstly remark that from the
theoretical point of view, it is challenging to propose a particular
choice of $C$ without having a specific cost function in mind,
which the thresholding approach inherently avoids. Therefore, one
possibility is
to use a large-scale simulation study to select a default value
of $C$ that works well across a range of signals.

With this in mind, we conducted the following simulation study.
For a given average number $N_{\mathrm{avg}} \in\{4, 8\}$ of change-points, we
simulated a Poisson number of change-points $N = \operatorname{Pois}(N_{\mathrm{avg}})$
and distributed them uniformly on $[0,1]$. At each change-point, we
introduced a jump whose height
had been drawn from the normal distribution with mean zero and
variance $\sigma^2_{\mathrm{jmp}} \in\{1, 3, 10\}$. We sampled the
thus-constructed function
at $T \in\{100, 200, 500, 1000, 2000\}$ equispaced points, and
contaminated it with Gaussian noise
with mean zero and variance one. Based on a large number of replicates, we
considered the quantity $|\hat{N} - N|$, where $\hat{N}$ was produced
by the WBS algorithm with threshold $\zeta_T = C \sqrt{2} \log^{1/2}
T$, and found
the value of $C$ that minimised it. The minimiser was sufficiently close
to $C = 1$ for us to use this value as the default one.

We add that our theoretical results do not permit a data-dependent
choice of the threshold constant $C$, so having a reliable default choice
is essential. The hope is that choosing such a default constant via
extensive simulation
should lead to good calibration of our method for a wide range of signals.

When the variance of $\varepsilon_t$ is unknown, we use $\zeta_T =
\hat{\sigma} C (2 \log T)^{1/2}$, where $\hat{\sigma}$ is the median
absolute deviation estimator of $\operatorname{Var}^{1/2}(\varepsilon_t)$.

Finally, we remark that in our comparative simulation study reported below,
we apply two threshold constants: the default value of $C = 1$ and
a higher value of $C = 1.3$. The latter is used for comparative purposes
as it was also used in the example considered in \citet{f14}.

Users with a preference for a method whose default parameters are not chosen
by simulation are encouraged to use the WBS method with the sSIC stopping
criterion described in Section~\ref{sec:ssic}, rather than with thresholding.
This method is also part of the simulation study below.

%s4.2 #&#
\subsection{Simulation study}
\label{sec:ss}

%t1 #&#
\begin{table}
\caption{Distribution of $\hat{N} - N$ for the various competing
methods and models, over 100 simulated sample paths. Also the average
mean-square error of the resulting estimate of $f_t$. Bold:
methods with the highest empirical frequency of $\hat{N} - N = 0$, and
those with frequencies within 10\% off the highest}
\label{tab1}
\begin{tabular*}{\textwidth}{@{\extracolsep{\fill}}lcd{3.0}d{2.0}d{2.0}d{2.0}d{2.0}d{2.0}d{2.0}c@{}}
\hline
& & \multicolumn{7}{c@{}}{$\bolds{\hat{N} - N}$} & \\[-6pt]
& & \multicolumn{7}{c@{}}{\hrulefill} & \\
\textbf{Method} & \textbf{Model} & \multicolumn{1}{c}{$\bolds{\le-3}$} & \multicolumn{1}{c}{$\bolds{-2}$} &
\multicolumn{1}{c}{$\bolds{-1}$} & \multicolumn{1}{c}{\textbf{0}} & \multicolumn{1}{c}{\textbf{1}} & \multicolumn{1}{c}{\textbf{2}} &
\multicolumn{1}{c}{$\bolds{\ge3}$} & \multicolumn{1}{c@{}}{\textbf{MSE}}\\
\hline
PELT &(1) & 0 & 0 & 0 & 8 & 9 & 9 & 74 & 4.3\\
B\&P & & 100 & 0 & 0 & 0 & 0 & 0 & 0 & 14.3\\
cumSeg & & 53 & 21 & 24 & 2 & 0 & 0 & 0 & 7.26\\
S3IB & & 0 & 5 & 42 & \textbf{51} & 1 & 1 & 0 & 2.55\\
SMUCE &  & 54 & 42 & 4 & 0 & 0 & 0 & 0 & 6.66\\
WBS $C = 1.0$ & & 0 & 0 & 24 & 38 & 15 & 18 & 5 & 2.77 \\
WBS $C = 1.3$ & & 1 & 13 & 78 & 8 & 0 & 0 & 0 & 3.02 \\
WBS sSIC & & 0 & 1 & 51 & \textbf{46} & 2 & 0 & 0 & 2.65 \\
BS $C = 1.0$ & & 0 & 1 & 40 & 39 & 16 & 2 & 2 & 3.12 \\
BS $C = 1.3$ & & 9 & 27 & 56 & 6 & 2 & 0 & 0 & 4.27 \\[3pt]
PELT & (2)& 0 & 0 & 0 & 15 & 11 & 20 & 54 & 79 $\times10^{-4}$\\
B\&P & & 99 & 1 & 0 & 0 & 0 & 0 & 0 & 399 $\times10^{-4}$\\
cumSeg & & 0 & 73 & 1 & 24 & 1 & 1 & 0 & 127 $\times10^{-4}$\\
S3IB & & 0 & 0 & 0 & \textbf{89} & 9 & 2 & 0 & 37 $\times10^{-4}$\\
SMUCE &  & 0 & 8 & 46 & 46 & 0 & 0 & 0 & 157 $\times10^{-4}$\\
WBS $C = 1.0$ & & 0 & 0 & 0 & 32 & 25 & 16 & 27 & 54 $\times10^{-4}$\\
WBS $C = 1.3$ & & 0 & 0 & 6 & \textbf{92} & 2 & 0 & 0 & 43 $\times
10^{-4}$\\
WBS sSIC & & 0 & 0 & 0 & \textbf{95} & 5 & 0 & 0 & 40 $\times10^{-4}$\\
BS $C = 1.0$ & & 0 & 0 & 30 & 49 & 16 & 4 & 1 & 75 $\times10^{-4}$\\
BS $C = 1.3$ & & 1 & 10 & 65 & 23 & 1 & 0 & 0 & 109 $\times10^{-4}$\\[3pt]
PELT &(3) & 0 & 0 & 3 & 11 & 16 & 17 & 53 & 2.08\\
B\&P & & 100 & 0 & 0 & 0 & 0 & 0 & 0 & 11.82\\
cumSeg & & 99 & 1 & 0 & 0 & 0 & 0 & 0 & 8.59\\
S3IB & & 34 & 34 & 18 & 14 & 0 & 0 & 0 & 1.96\\
SMUCE &  & 63 & 28 & 8 & 1 & 0 & 0 & 0 & 4.35\\
WBS $C = 1.0$ & & 0 & 9 & 22 & \textbf{32} & 21 & 13 & 3 & 1.67\\
WBS $C = 1.3$ & & 15 & 41 & 32 & 12 & 0 & 0 & 0 & 1.91\\
WBS sSIC & & 7 & 28 & 23 & \textbf{33} & 6 & 1 & 2 & 1.62\\
BS $C = 1.0$ & & 10 & 30 & 26 & 19 & 13 & 2 & 0 & 2.34\\
BS $C = 1.3$ & & 80 & 18 & 2 & 0 & 0 & 0 & 0 & 3.99\\[3pt]
PELT &(4) & 0 & 0 & 0 & 38 & 28 & 18 & 16 & $55 \times10^{-3}$\\
B\&P & & 100 & 0 & 0 & 0 & 0 & 0 & 0 & $251 \times10^{-3}$\\
cumSeg & & 100 & 0 & 0 & 0 & 0 & 0 & 0 & $251 \times10^{-3}$\\
S3IB & & 36 & 16 & 1 & 47 & 0 & 0 & 0 & $116 \times10^{-3}$\\
SMUCE &  & 98 & 1 & 0 & 1 & 0 & 0 & 0 & $215 \times10^{-3}$\\
WBS $C = 1.0$ & & 0 & 1 & 7 & \textbf{77} & 11 & 2 & 2 & $51 \times
10^{-3}$\\
WBS $C = 1.3$ & & 22 & 11 & 28 & 38 & 1 & 0 & 0 & $80 \times10^{-3}$\\
WBS sSIC & & 4 & 1 & 4 & \textbf{80} & 7 & 4 & 0 & $55 \times10^{-3}$\\
BS $C = 1.0$ & & 49 & 9 & 19 & 14 & 9 & 0 & 0 & $129 \times10^{-3}$\\
BS $C = 1.3$ & & 94 & 3 & 3 & 0 & 0 & 0 & 0 & $210 \times10^{-3}$\\
\hline
\end{tabular*}
\end{table}
\setcounter{table}{0}
\begin{table}
\caption{(Continued)}
\begin{tabular*}{\textwidth}{@{\extracolsep{\fill}}lcd{3.0}d{2.0}d{2.0}d{2.0}d{2.0}d{2.0}d{2.0}c@{}}
\hline
& & \multicolumn{7}{c@{}}{$\bolds{\hat{N} - N}$} & \\[-6pt]
& & \multicolumn{7}{c@{}}{\hrulefill} & \\
\textbf{Method} & \textbf{Model} & \multicolumn{1}{c}{$\bolds{\le-3}$} & \multicolumn{1}{c}{$\bolds{-2}$} &
\multicolumn{1}{c}{$\bolds{-1}$} & \multicolumn{1}{c}{\textbf{0}} & \multicolumn{1}{c}{\textbf{1}} & \multicolumn{1}{c}{\textbf{2}} &
\multicolumn{1}{c}{$\bolds{\ge3}$} & \multicolumn{1}{c@{}}{\textbf{MSE}}\\
\hline
PELT & (5)& 0 & 0 & 0 & 34 & 24 & 19 & 23 & $26 \times10^{-3}$\\
B\&P & & 100 & 0 & 0 & 0 & 0 & 0 & 0 & $554 \times10^{-3}$\\
cumSeg & & 3 & 1 & 11 & 77 & 8 & 0 & 0 & $63 \times10^{-3}$\\
S3IB & & 97 & 1 & 2 & 0 & 0 & 0 & 0 & $210 \times10^{-3}$\\
SMUCE &  & 64 & 17 & 11 & 8 & 0 & 0 & 0 & $185 \times10^{-3}$\\
WBS $C = 1.0$ & & 0 & 0 & 0 & 63 & 31 & 4 & 2 & $24 \times10^{-3}$\\
WBS $C = 1.3$ & & 0 & 0 & 4 & \textbf{87} & 9 & 0 & 0 & $27 \times
10^{-3}$\\
WBS sSIC & & 0 & 0 & 0 & 61 & 35 & 4 & 0 & $23 \times10^{-3}$\\
BS $C = 1.0$ & & 0 & 0 & 0 & 79 & 20 & 1 & 0 & $24 \times10^{-3}$\\
BS $C = 1.3$ & & 0 & 0 & 5 & \textbf{88} & 7 & 0 & 0 & $27 \times10^{-3}$\\
\hline
\end{tabular*}
\end{table}

%t2 #&#
\begin{table}
\caption{Summary statistics for the empirical distribution of $\hat{N}
- N$
for the various competing methods and models, over 100 simulated sample
paths}
\label{tab2}
\begin{tabular*}{\textwidth}{@{\extracolsep{\fill}}lcd{3.0}d{3.2}d{3.1}d{3.2}d{2.2}d{3.0}@{}}
\hline
& & \multicolumn{6}{c@{}}{\textbf{Summary of} $\bolds{\hat{N} - N}$} \\[-6pt]
& & \multicolumn{6}{c@{}}{\hrulefill}  \\
\textbf{Method} & \textbf{Model} & \multicolumn{1}{c}{\textbf{Min.}} & \multicolumn{1}{c}{\textbf{1st Qu.}} &
\multicolumn{1}{c}{\textbf{Median}} & \multicolumn{1}{c}{\textbf{Mean}} & \multicolumn{1}{c}{\textbf{3rd Qu.}} & \multicolumn{1}{c@{}}{\textbf{Max.}} \\
\hline
PELT &(1) & 0 & 2 & 4 & 4.78 & 7 & 18 \\
B\&P & & -7 & -7 & -7 & -7 & -7 & -7 \\
cumSeg & & -7 & -4 & -3 & -2.84 & -1 & 0 \\
S3IB & & -2 & -1 & 0 & -0.49 & 0 & 2 \\
SMUCE &  & -4 & -3 & -3 & -2.63 & -2 & -1\\
WBS $C = 1.0$ & & -1 & 0 & 0 & 0.45 & 1 & 5 \\
WBS $C = 1.3$ & & -3 & -1 & -1 & -1.07 & -1 & 0 \\
WBS sSIC & & -2 & -1 & -1 & -0.51 & 0 & 1 \\
BS $C = 1.0$ & & -2 & -1 & 0 & -0.16 & 0 & 3 \\
BS $C = 1.3$ & & -3 & -2 & -1 & -1.35 & -1 & 1 \\[3pt]
PELT & (2)& 0 & 1 & 3 & 3.39 & 5 & 10 \\
B\&P & & -3 & -3 & -3 & -2.99 & -3 & -2 \\
cumSeg & & -2 & -2 & -2 & -1.44 & 0 & 2 \\
S3IB & & 0 & 0 & 0 & 0.13 & 0 & 2 \\
SMUCE &  & -2 & -1 & -1 & -0.62 & 0 & 0\\
WBS $C = 1.0$ & & 0 & 0 & 1 & 1.56 & 3 & 5 \\
WBS $C = 1.3$ & & -1 & 0 & 0 & -0.04 & 0 & 1 \\
WBS sSIC & & 0 & 0 & 0 & 0.05 & 0 & 1 \\
BS $C = 1.0$ & & -1 & -1 & 0 & -0.03 & 0 & 3 \\
BS $C = 1.3$ & & -3 & -1 & -1 & -0.87 & -1 & 1 \\[3pt]
PELT &(3) & -1 & 1 & 3 & 3.17 & 4.25 & 12 \\
B\&P & & -13 & -13 & -12.5 & -12.44 & -12 & -10 \\
cumSeg & & -13 & -13 & -9 & -9 & -5 & -2 \\
S3IB & & -12 & -3 & -2 & -2.15 & -1 & 0 \\
SMUCE &  & -6 & -4 & -3 & -2.95 & -2 & 0\\
WBS $C = 1.0$ & & -2 & -1 & 0 & 0.16 & 1 & 3 \\
WBS $C = 1.3$ & & -5 & -2 & -2 & -1.64 & -1 & 0 \\
WBS sSIC & & -5 & -2 & -1 & -0.88 & 0 & 4 \\
BS $C = 1.0$ & & -5 & -2 & -1 & -1.03 & 0 & 2 \\
BS $C = 1.3$ & & -8 & -4 & -3 & -3.56 & -3 & -1 \\[3pt]
PELT &(4) & 0 & 0 & 1 & 1.26 & 2 & 7 \\
B\&P & & -13 & -13 & -13 & -12.98 & -13 & -12 \\
cumSeg & & -13 & -13 & -13 & -12.94 & -13 & -11 \\
S3IB & & -13 & -11.25 & -2 & -4.17 & 0 & 0 \\
SMUCE &  & -12 & -8 & -6 & -6.35 & -5 & 0\\
WBS $C = 1.0$ & & -2 & 0 & 0 & 0.12 & 0 & 3 \\
WBS $C = 1.3$ & & -8 & -2 & -1 & -1.46 & 0 & 1 \\
WBS sSIC & & -12 & 0 & 0 & -0.21 & 0 & 2 \\
BS $C = 1.0$ & & -11 & -4.25 & -2 & -2.70 & -1 & 1 \\
BS $C = 1.3$ & & -13 & -11 & -9 & -8.42 & -7 & -1 \\
\hline
\end{tabular*}
\end{table}
\setcounter{table}{1}
\begin{table}
\caption{(Continued)}
\begin{tabular*}{\textwidth}{@{\extracolsep{\fill}}lcd{3.0}d{3.2}d{3.1}d{3.2}d{2.2}d{3.0}@{}}
\hline
& & \multicolumn{6}{c@{}}{\textbf{Summary of} $\bolds{\hat{N} - N}$} \\[-6pt]
& & \multicolumn{6}{c@{}}{\hrulefill}  \\
\textbf{Method} & \textbf{Model} & \multicolumn{1}{c}{\textbf{Min.}} & \multicolumn{1}{c}{\textbf{1st Qu.}} &
\multicolumn{1}{c}{\textbf{Median}} & \multicolumn{1}{c}{\textbf{Mean}} & \multicolumn{1}{c}{\textbf{3rd Qu.}} & \multicolumn{1}{c@{}}{\textbf{Max.}} \\
\hline
PELT &(5) & 0 & 0 & 1 & 1.55 & 2 & 8 \\
B\&P & & -9 & -9 & -9 & -9 & -9 & -9 \\
cumSeg & & -5 & 0 & 0 & -0.17 & 0 & 1 \\
S3IB & & -7 & -7 & -6 & -5.71 & -5 & -1 \\
SMUCE &  & -6 & -4 & -3 & -2.85 & -2 & 0\\
WBS $C = 1.0$ & & 0 & 0 & 0 & 0.46 & 1 & 4 \\
WBS $C = 1.3$ & & -1 & 0 & 0 & 0.05 & 0 & 1 \\
WBS sSIC & & 0 & 0 & 0 & 0.43 & 1 & 2 \\
BS $C = 1.0$ & & 0 & 0 & 0 & 0.22 & 0 & 2 \\
BS $C = 1.3$ & & -1 & 0 & 0 & 0.02 & 0 & 1 \\
\hline
\end{tabular*}
\end{table}

In this section, we compare the performance of WBS (and BS)
against the best available competitors implemented in R packages,
most of which are publicly available on CRAN. The competing packages are:
\texttt{strucchange}, which implements the multiple change-point
detection method of \citet{bp03}, \texttt{Segmentor3IsBack},
which implements the method of \citet{r10} with the model selection
methodology from \citet{l05}, \texttt{changepoint}, which implements the
PELT methodology of \citet{kfe12}, \texttt{cumSeg}, which implements
the method from \citet{ma11}, and \texttt{stepR}, which implements the SMUCE
method of \citet{fms14}.
In the remainder of this section,
we refer to these methods as, respectively, B\&P, S3IB, PELT, cumSeg, and
SMUCE. Appendix~\ref{sec:simextra} provides an extra discussion of how
these methods were used in our simulation study.
With the exception of \texttt{stepR}, which
is available from \href{http://www.stochastik.math.uni-goettingen.de/smuce}{http://www.stochastik.math.uni-goettingen.de/\break smuce}
at the time of writing, the remaining packages are available on CRAN.

In this section, the WBS algorithm uses the default value of $M = 5000$ random
draws. In the thresholding stopping rule, we use the threshold
$\zeta_T = C  \hat{\sigma} \sqrt{2 \log T}$, where $\hat{\sigma}$ is the
median absolute deviation estimator of $\sigma$ suitable for i.i.d.
Gaussian noise,
$T$ is the sample size, and the constant $C$ is set to 1 and $1.3$ as motivated
earlier. The WBS method combined with the sSIC stopping criterion is
referred to
as `WBS sSIC' and uses $\alpha= 1.01$, again as justified earlier, and
$K = 20$.
The BS method uses the same thresholds as WBS, for comparability.

Our test signals, fully specified in Appendix~\ref{sec:simextra}
along with the sample sizes and noise standard deviations used,
are (1) \texttt{blocks}, (2) \texttt{fms},
(3) \texttt{mix},
(4) \texttt{teeth10},
and (5) \texttt{stairs10}. Tables~\ref{tab1} and \ref{tab2} show the results.
We describe the performance of each method below.

\begin{longlist}[\quad]
\item[B\&P.] The B\&P method performs poorly, which may be partly
due to the default minimum segment size set to 15\% of the sample
size, an assumption violated by several of our test signals. However,
resetting this parameter to 1 or even 1\% of the sample size resulted in
exceptionally slow computation times, which prevented us from reporting the
results in our comparisons.

\item[S3IB.] This method offers excellent performance for the \texttt{blocks}
signal, and very good performance for the \texttt{fms} signal. The
\texttt{mix}
signal is more challenging, and the S3IB method does not perform well here,
with a tendency to underestimate the number of change-points, sometimes by
as many as 12. Performance is rather average for the \texttt{teeth10} signal,
and systematically poor for the \texttt{stairs10} signal.

\item[PELT.] The PELT method has a tendency to overestimate the number of
change-points, which is apparent in all of the examples studied.

\item[cumSeg.] Apart from the \texttt{stairs10} signal
for which it offers acceptable performance, the cumSeg method tends to heavily
underestimate the number of change-points.

\item[SMUCE.] The SMUCE method tends to underestimate the true number
of change-points. However, its performance for the \texttt{fms} signal
is acceptable.

\item[BS.] For $C=1$, the method performs acceptably for the \texttt{blocks} and \texttt{stairs10} signals, has rather average performance for the
\texttt{fms} and \texttt{mix} signals, and performs poorly for \texttt{teeth10}.
For $C=1.3$, performance is excellent for the \texttt{stairs10} signal;
otherwise poor. Overall, our test signals clearly demonstrate the
shortcomings of classical binary segmentation.

\item[WBS.]
With the threshold constant $C = 1$, WBS works well for the \texttt{blocks} and \texttt{stairs10} signals, although in both cases it is behind the
best performers. For the \texttt{fms} signal, it tends to overestimate the
number of change-points, although not by many. It offers (relatively)
excellent performance for \texttt{mix} and \texttt{teeth10}.

For $C = 1.3$, WBS performs excellently for \texttt{fms} and \texttt{stairs10},
while it underestimates the number of change-points for the other signals,
although again, not by many.

WBS sSIC performs the best or very close to the best for all signals
bar \texttt{stairs10}; however,
for the latter, if it overestimates the number of change-points, then
it does so mostly
by one change-point only. If one overall `winner' were to be chosen out
of the
methods studied, it would clearly have to be WBS sSIC.
\end{longlist}

Our overall recommendation is to use WBS sSIC first. If the visual inspection
of the residuals from the fit reveals any obvious patterns neglected by
WBS sSIC, then WBS with $C = 1.3$ should be used next. Since the latter
has a tendency
to underestimate the number of change-points, the hope is that it does not
detect any spurious ones. If patterns in residuals remain, WBS with
$C=1$ should be
used next.

Furthermore, Appendix~\ref{sec:linear} contains a small-scale simulation
study and brief discussion regarding the performance of WBS in the presence
of linear trends.

%s5 #&#
\section{Real data example}
\label{sec:real}

In this section, we apply the WBS method to the detection of trends
in the S\&P 500 index. We consider the time series of log-returns on
the daily closing values of S\&P 500, of length $T = 2000$ (i.e., approximately
8 trading years) ending 26 October 2012. We then remove the volatility
of this series by fitting the $\operatorname{GARCH}(1,1)$ model with Gaussian innovations,
and apply the WBS procedure to the residuals $X_t$ from the fit, both
with the
thresholding and the sSIC stopping criteria. To obtain a more complete picture
of the estimated change-point structure, it is instructive to carry out
the WBS
procedure for a range of thresholds~$\zeta_T$.

\begin{figure}

\includegraphics{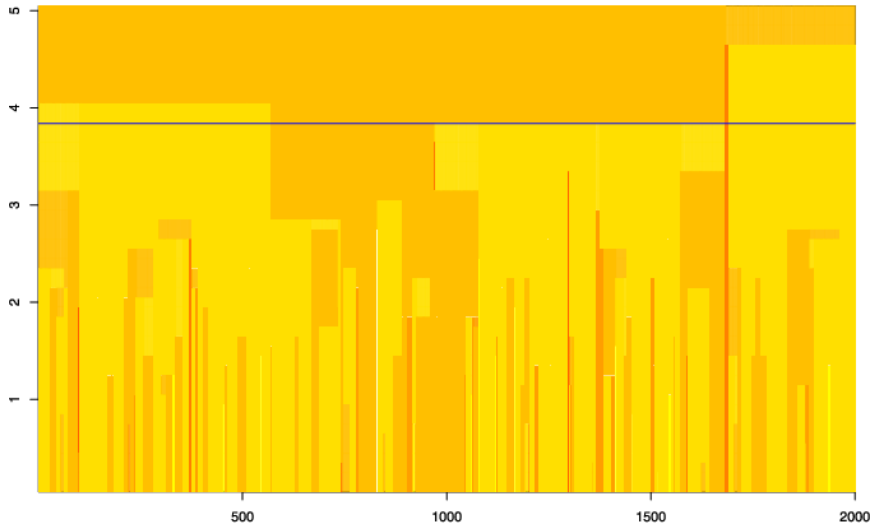}

\caption{Time-threshold map of the WBS acting on the series $X_t$ from
Section \protect\ref{sec:real}. The horizontal blue line is the
threshold $\zeta_T \approx3.83$.}\label{fig:map}
\end{figure}

The results, for $\zeta_T$ changing from 0 to 5, are presented in the
`time-threshold map' [see \citet{f12} for more details of this generic
concept] in Figure~\ref{fig:map}. The map should be read as follows.
The $x$-coordinates of
the vertical lines indicate the estimated change-point locations,
detected for the range of
thresholds equal to the range of the given line on the $y$-axis. For example,
for $\zeta_T = \hat{\sigma}(2\log T)^{1/2} \approx3.83$, we
have 5 estimated change-points, since the horizontal blue line
(corresponding to
$\zeta_T = 3.83$) in Figure~\ref{fig:map} crosses 5 vertical lines. The
5 estimated
change-points are concentrated in or around 3 separate locations.

%f3 #&#

%f4 #&#
\begin{figure}

\includegraphics{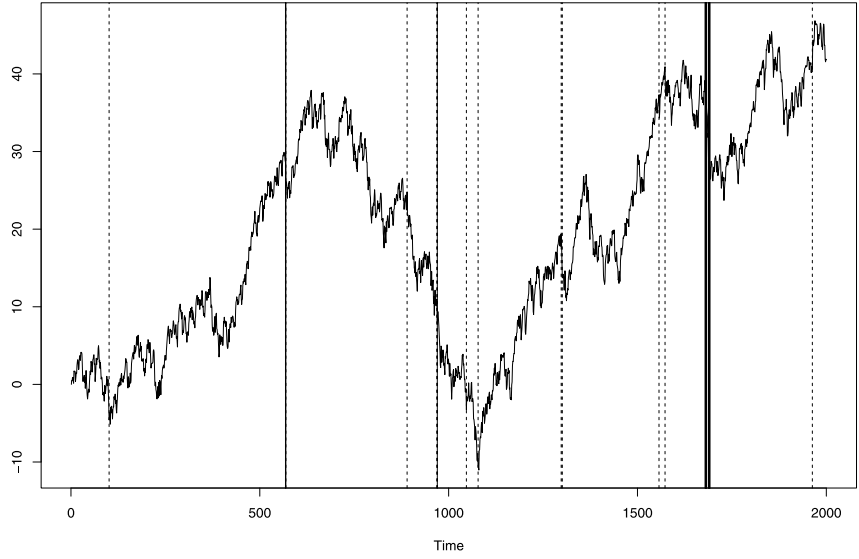}

\caption{Cumulative sum of $X_t$, change-points corresponding to
sSIC (thick solid vertical lines), $\zeta_T = 3.83$
(thin and thick solid vertical lines), $\zeta_T = 3.1$
(all vertical lines).}\label{fig:31}
\end{figure}

Figure~\ref{fig:31} shows the corresponding cumulative sum of the
residuals from
the GARCH fit (which can be viewed as the logged S\&P 500 index with its
volatility removed), with the estimated change-point locations
corresponding to the
thresholds $\zeta_T = 3.83$ and $\zeta_T = 3.1$, as well as the sSIC criterion.
Interestingly, the sSIC criterion estimates only 2 change-points, both
concentrated
around time $t = 1700$.

As with any other financial data, it is difficult to speak
of the number of estimated change-points being `right' or `wrong' here:
for example,
some more frequent traders may naturally be more interested in trend
changes on the scale of
weeks or months, rather than years, in which case a lower threshold
might be more
suitable. However, it is interesting to observe that both the sSIC
criterion, the most
accurate estimator of $\hat{N}$ from our simulation study, and the
thresholding criterion with
$\zeta_T = 3.83$, which corresponds to the threshold constant $C = 1$
and tended to slightly overestimate
the number of change-points in the simulation study, point to a rather
low number of
estimated change-points in this example.

\begin{appendix}\label{app}
%s6 #&#
\section{Proofs}

\begin{pf*}{Proof of Theorem~\ref{th:sbs}}
We first introduce some notation. Denoting $n = e-s+1$, we define
%
%e5 #&#
%e6 #&#
\begin{eqnarray}
\label{eq:xtilde} \tilde{X}_{s,e}^b & = & \sqrt{
\frac{e-b}{n(b-s+1)}}\sum_{t=s}^b
X_t - \sqrt{\frac{b-s+1}{n(e-b)}}\sum_{t=b+1}^e
X_t,
\\
\label{eq:ftilde} \tilde{f}_{s,e}^b & = & \sqrt{
\frac{e-b}{n(b-s+1)}}\sum_{t=s}^b
f_t - \sqrt{\frac{b-s+1}{n(e-b)}}\sum_{t=b+1}^e
f_t.
\end{eqnarray}
Let $s, e$ satisfy
%
%e7 #&#
\begin{equation}
\label{eq:cond0} \eta_{p_0} \le s < \eta_{p_0+1} < \cdots<
\eta_{p_0+q} < e \le\eta_{p_0+q+1}
\end{equation}
for $0 \le p_0 \le N-q$, which will be the case at all stages of the
algorithm while there are
still undetected change-points remaining. In Lemmas \ref{lem:near}--\ref
{lem:howlarge}, we impose the following conditions:
%
%e8 #&#
%e9 #&#
\begin{eqnarray}
\label{eq:cond1}&& s < \eta_{p_0+r} - C\delta_T <
\eta_{p_0+r} + C \delta_T < e \qquad\mbox {for some } 1 \le r \le q,
\\
\label{eq:cond2}&& \max\bigl(\min(\eta_{p_0+1}-s, s-\eta_{p_0}),
\min(\eta_{p_0+q+1}-e, e-\eta _{p_0+q})\bigr) \le C
\epsilon_T.
\end{eqnarray}
Both (\ref{eq:cond1}) and (\ref{eq:cond2}) hold throughout the
algorithm for all those segments
starting at $s$ and ending at $e$ which contain previously undetected
change-points. As Lemma~\ref{lem:howsmall} concerns the case where all change-points have been
detected, it does not use either of these
conditions.

We also introduce a set $A_T$ defined by
%
%e10 #&#
\begin{equation}
\label{eq:defat} A_T = \Biggl\{ \Biggl\llvert (e-b+1)^{-1/2}
\sum_{i=b}^e \varepsilon_i
\Biggr\rrvert < \lambda_2\ \forall 1 \le b \le e \le T \Biggr\}.
\end{equation}
Note that by Bonferroni's inequality, $P(A_T) \ge1 - C T^{-1}$ for
$\lambda_2 \ge(6 \log T)^{1/2}$, where $C$ is
a positive constant.

Before presenting the formal proof, we informally discuss some of its
aspects to facilitate
understanding.

\textit{Informal discussion of some aspects of the proof.}
The performance of the binary segmentation
algorithm analysed in Theorem~\ref{th:sbs} can be seen as
`deterministic on
a random set whose probability approaches one',
in the sense that for a $T$ large enough and in a certain subset
of the probability space whose probability approaches one,
the algorithm is guaranteed to detect all true change-points
before being stopped at the right time by the application of threshold
$\zeta_T$. We further clarify this observation below.

Heuristically speaking, on the set $A_T \cap B_T$, where $A_T$ is
defined in (\ref{eq:defat})
and $B_T$ in Lemma~\ref{lem:bonf}, the innovations $\varepsilon_t$
are well behaved in the sense that the empirical CUSUM statistics
$\tilde{X}_{s,e}^b$ are uniformly close to the\vspace*{-1pt} corresponding unobserved true
quantities $\tilde{f}_{s,e}^b$ in the particular sense described
in Lemmas \ref{lem:bonf} and~\ref{lem:vnear}. It is this closeness that
causes the following
behaviour: if there are still previously undetected change-points within
the current interval $[s,e]$ (by which we mean that there are change-points
for which there is no estimated change-point within the distance of $C
\epsilon_T$),
and $[s,e]$ satisfies (\ref{eq:cond1}) and (\ref{eq:cond2}),
then (i) by Lemma~\ref{lem:vnear}, $b_0 = \arg\max_{t: s\le t<e}|\tilde
{X}_{s,e}^t|$
falls within the distance of $C \epsilon_T$ of one of the previously undetected
change-points in $[s,e]$ (denote that change-point here by $\eta
_{p_0+r}$), and~(ii)
by Lemma~\ref{lem:howlarge}, we have $|\tilde{X}_{s,e}^{b_0}| > \zeta_T$.

The consequence of (i) and (ii) is that $b_0$ passes the thresholding
test for the
significance of a change-point and is from now on considered to be an
estimate of $\eta_{p_0+r}$. Note that the assignment of $b_0$ to $\eta
_{p_0+r}$ is unambiguous:
$b_0$ cannot be an estimate of any of the other change-points as they
are too far; the nearest left- or right-neighbour of $\eta_{p_0+r}$
is at a distance of no less than $\delta_T$ of it, which means not
nearer than $\delta_T - C \epsilon_T$ from $b_0$, which is orders
of magnitude larger than $C \epsilon_T$ as specified in Theorem~\ref{th:sbs}
and its assumptions.

As a consequence, the procedure then moves on to operate on the intervals
$[s,b_0]$ and $[b_0,e]$. Without loss of generality, suppose there are
previously undetected change-points on $[s,b_0]$. We now demonstrate
that (\ref{eq:cond1}) and (\ref{eq:cond2}) hold for that interval.
Since $b_0$ is close
to $\eta_{p_0+r}$ (which is `previously detected'), it must be
far from all other true change-points in the sense described
in the previous paragraph. In particular, for any previously undetected
change-point $\eta_{p_0+r'} \in[s,b_0]$, we must have $b_0 - \eta
_{p_0+r'} \ge
\delta_T - C \epsilon_T$ by the argument from the previous paragraph,
which is larger than $C \delta_T$ for some $C > 0$, by the assumptions
of Theorem~\ref{th:sbs}. Hence $[s,b_0]$ satisfies (\ref{eq:cond1}).

Similarly, $[s,b_0]$ satisfies (\ref{eq:cond2}) as $b_0$ is within the
distance of
$C\epsilon_T$ of one of its neighbouring change-points, namely $\eta_{p_0+r}$.

Thus (\ref{eq:cond1}) and (\ref{eq:cond2}) are both valid as the
algorithm progresses for any interval
$[s,e]$ on which there are still previously undetected change-points. Therefore,
for a large enough $T$ and on $A_T \cap B_T$, all change-points will be
detected one by one. At that point, by Lemma~\ref{lem:howsmall}, the statistics
$|\tilde{X}_{s,e}^b|$ will become uniformly smaller than the threshold
$\zeta_T$ and the algorithm will stop.

We are now in a position to turn to the formal proof, which is split
into a number
of lemmas.

%le6.1 #&#
\begin{lemma}
\label{lem:bonf}
Let $X_t$ follow model (\ref{eq:mod}), and let the assumptions of
Theorem~\ref{th:sbs}
hold. Let $\tilde{X}_{s,e}^b$ and $\tilde{f}_{s,e}^b$ be defined by
(\ref{eq:xtilde}) and (\ref{eq:ftilde}), respectively. We then have
$P(B_T) \ge1 - C T^{-1}$, where
\[
B_T = \Bigl\{ \max_{s,b,e  \dvtx 1 \le s \le b < e \le T}\bigl |\tilde{X}_{s,e}^b
- \tilde{f}_{s,e}^b\bigr| \le\lambda_1 \Bigr\},
\]
$\lambda_1 \ge\sqrt{8 \log T}$, and $C$ is a positive constant.
\end{lemma}

\begin{pf}The proof proceeds via a simple Bonferroni inequality,
\[
1 - P(B_T) \le\sum_{s,b,e} P\bigl(|Z| >
\lambda_1\bigr) \le T^3 \frac{\phi
_Z(\lambda_1)}{\lambda_1} \le
\frac{C}{T},
\]
where $Z$ is a standard normal and $\phi_Z(\cdot)$ is its p.d.f.
\end{pf}

We conjecture that more accurate bounds for $\lambda_1$ of Lemma~\ref
{lem:bonf} and
$\lambda_2$ in formula (\ref{eq:defat}) can be obtained, for example,
using techniques as in
\citet{twg07}, \citet{aj13}, or especially Lemma~1 of \citet{y88}.
However, we note that
even with the use of the suboptimal Bonferroni inequality, $\lambda_1$
and $\lambda_2$
are already rate-optimal, which is what matters for the rates of
convergence in
Theorems \ref{th:sbs}--\ref{th:ssic}. Improving the multiplicative
constants in $\lambda_1$ and $\lambda_2$
would bring no further practical benefits in terms of choosing the
stopping criterion
for BS or WBS, the main reason for this being that the result of Lemma~\ref{lem:howsmall}
(and its equivalent in the proof of Theorem~\ref{th:wbs}) is dependent
on a different
constant $C$ anyway, which is not straightforward to evaluate in theory.

%le6.2 #&#
\begin{lemma}
\label{lem:near}
Let $X_t$ follow model (\ref{eq:mod}) and let the assumptions of
Theorem~\ref{th:sbs}
hold. Let $\tilde{X}_{s,e}^b$ and $\tilde{f}_{s,e}^b$ be defined by
(\ref{eq:xtilde}) and (\ref{eq:ftilde}), respectively.
Assume (\ref{eq:cond0}), (\ref{eq:cond1}), and~(\ref{eq:cond2}).
On set $B_T$ of Lemma~\ref{lem:bonf}, the following
holds. For $b = \arg\max_{t: s \le t < e}|\tilde{X}_{s,e}^t|$, there
exists $1 \le r \le q$ such that for large $T$, $|b - \eta_{p_0+r}| \le
C_1\gamma_T$ with
$\gamma_T = T^{1/2} \lambda_1 / f'_{p_0+r}$
($\lambda_1$ as in Lemma~\ref{lem:bonf}). In addition, $|\tilde
{f}_{s,e}^t|$ must then have a local
maximum at $t = \eta_{p_0+r}$, and we must have
\[
\frac{|\tilde{f}_{s,e}^{\eta_{p_0+r}}|}{\max_{t \dvtx s\le t<e}|\tilde
{f}_{s,e}^t|} \ge C_2,
\]
where $C_1, C_2$ are positive constants.
\end{lemma}

\begin{pf}
We first note that $\gamma_T = o(\delta_T)$ since
$1/2 + \varpi< 1/2 + 2\Theta- 3/2 \le\Theta$. Note also
that $\delta_T T^{-1/2}\underline{f}_T \ge C T^\varphi$ for $C, \varphi
$ positive.
Let $b_1 =\break  \arg\max_{t: s\le t<e} |\tilde{f}_{s,e}^t|$. From Lemma~\ref
{lem:bonf},
we have
%
%e11 #&#
\begin{equation}
\label{eq:contr} \bigl|\tilde{f}_{s,e}^{b_1}\bigr| \le\bigl|
\tilde{X}_{s,e}^{b_1}\bigr| + \lambda_1 \le \bigl|
\tilde{X}_{s,e}^{b}\bigr| + \lambda_1 \le\bigl|
\tilde{f}_{s,e}^{b}\bigr| + 2\lambda_1.
\end{equation}
Assume $b \in(\eta_{p_0+r} + C\gamma_T, \eta_{p_0+r+1} - C \gamma_T)$
for a large enough constant $C$
for some $r$ and
w.l.o.g. $\tilde{f}_{s,e}^b > 0$. From Lemma~2.2 in Venkatraman (\citeyear{venkatraman1993}),
$\tilde{f}_{s,e}^t$
is either monotonic or decreasing and then increasing on $[\eta
_{p_0+r}, \eta_{p_0+r+1}]$
and $\max(\tilde{f}_{s,e}^{\eta_{p_0+r}}, \tilde{f}_{s,e}^{\eta
_{p_0+r+1}}) > \tilde{f}_{s,e}^b$.
If $\tilde{f}_{s,e}^b$ locally decreases at $b$, then $\tilde
{f}_{s,e}^{\eta_{p_0+r}} >
\tilde{f}_{s,e}^b$, and arguing exactly as in Lemma~2 of \citet{cf12},
there exists
$b' \in(\eta_{p_0+r}, \eta_{p_0+r} + C\gamma_T]$ such that $\tilde
{f}_{s,e}^{\eta_{p_0+r}} \ge
\tilde{f}_{s,e}^{b'} + 2\lambda_1$. This would in turn lead to $|\tilde
{f}_{s,e}^{b_1}| >
|\tilde{f}_{s,e}^{b}| + 2\lambda_1$, a contradiction of (\ref
{eq:contr}). Similar arguments apply
if $\tilde{f}_{s,e}^b$ locally increases at $b$.

Let $r$ be as in the statement of this lemma. Then $|\tilde
{f}_{s,e}^{\eta_{p_0+r}}|$ must be a
local maximum, as if it were not, we would have $\max(|\tilde
{f}_{s,e}^{\eta_{p_0+r-1}}|,
|\tilde{f}_{s,e}^{\eta_{p_0+r+1}}|) > |\tilde{f}_{s,e}^{\eta
_{p_0+r}}|$, and arguing exactly as
above, this maximum would have to be sufficiently larger than $|\tilde
{f}_{s,e}^{\eta_{p_0+r}}|$
for $b$ to fall near the change-point achieving this maximum, rather
than near $\eta_{p_0+r}$,
which is a contradiction.

Finally, using the same argumentation again, $|\tilde{f}_{s,e}^{\eta
_{p_0+r}}| / \max_{t: s\le t<e}|\tilde{f}_{s,e}^t|$
must be bounded from below, as if were not, then recalling that\break $\max_{t: s\le t<e}|
\tilde{f}_{s,e}^t| \ge C \delta_T T^{-1/2} \underline{f}_T$
by Lemma~1 of \citet{cf12}, $b$ would have to fall near the change-point
achieving this maximum, rather than
near $\eta_{p_0+r}$, which is again a contradiction. This completes the
proof of the lemma.
\end{pf}

%le6.3 #&#
\begin{lemma}
\label{lem:vnear}
Let the conditions of Lemma~\ref{lem:near} hold, and let the notation
be as in that lemma.
On set $B_T \cap A_T$, where $B_T$ is defined in Lemma~\ref{lem:bonf}
and $A_T$ in (\ref{eq:defat}),
we have for large $T$, $|b - \eta_{p_0+r}| \le C \epsilon_T$, where
$\epsilon_T =
\lambda_2^2 T^2 \times\break \delta_T^{-2} (f'_{p_0+r})^{-2}$
and $C$ is a positive constant.
\end{lemma}

\begin{pf} Let $\langle\cdot, \cdot\rangle$ denote the inner product
between two vectors. Let $\psi_{s,e}^d$ be a vector whose elements $\psi
_{s,e,t}^d$ are constant
and positive for $t = s, \ldots, d$, constant and negative for $t =
d+1, \ldots, e$, sum to zero
and such that their squares sum to one. Then it is easy to see that
$\tilde{X}_{s,e}^d = \sum_{t=s}^e \psi_{s,e,t}^d X_t = \langle\psi
_{s,e}^d, X \rangle$
and similarly $\tilde{f}_{s,e}^d = \langle\psi_{s,e}^d, f \rangle$.
For any vector $v$
supported on $[s,e]$, we have $\arg\max_{d:s\le d<e} |\langle\psi
_{s,e}^d, v \rangle| = \arg\min_{d:s\le d<e}
\sum_{t=s}^e (v_t - \bar{v}_{s,e,t}^d)^2$, where $\bar{v}_{s,e}^d$ is
an orthogonal
projection of $v$ on the space of step functions constant on $s, \ldots, d$ and constant
on $d+1, \ldots, e$; this is immediate by noting that $\bar{v}_{s,e}^d
= \bar{v} + \langle v, \psi_{s,e}^d\rangle
\psi_{s,e}^d$, where $\bar{v}$ is the mean of $v$. From Lemma~\ref{lem:near},
\[
\sum_{t=s}^e \bigl(X_t -
\bar{X}_{s,e,t}^b\bigr)^2 \le\sum
_{t=s}^e \bigl(X_t - \bar
{f}_{s,e,t}^{\eta_{p_0+r}}\bigr)^2.
\]
Therefore, if it can be shown that for a certain $\epsilon_T < C_1
\gamma_T$, we have
%
%e12 #&#
\begin{equation}
\label{eq:away} \sum_{t=s}^e
\bigl(X_t - \bar{X}_{s,e,t}^d\bigr)^2
> \sum_{t=s}^e \bigl(X_t -
\bar {f}_{s,e,t}^{\eta_{p_0+r}}\bigr)^2
\end{equation}
as long as
%
%e13 #&#
\begin{equation}
\label{eq:restr} \epsilon_T < |d - \eta_{p_0+r}| \le
C_1 \gamma_T,
\end{equation}
then this would prove that necessarily,
$|b - \eta_{p_0+r}| \le\epsilon_T$. Recalling that $X_t = f_t +
\varepsilon_t$, (\ref{eq:away}) is
equivalent to
\[
2 \sum_{t=s}^e \varepsilon_t
\bigl(\bar{X}_{s,e,t}^d - \bar{f}_{s,e,t}^{\eta
_{p_0+r}}
\bigr) < \sum_{t=s}^e \bigl(f_t
- \bar{X}_{s,e,t}^d\bigr)^2 - \sum
_{t=s}^e \bigl(f_t - \bar
{f}_{s,e,t}^{\eta_{p_0+r}}\bigr)^2,
\]
and implied by
%
%e14 #&#
\begin{equation}
\label{eq:stochdet} 2 \sum_{t=s}^e
\varepsilon_t \bigl(\bar{X}_{s,e,t}^d -
\bar{f}_{s,e,t}^{\eta
_{p_0+r}}\bigr) < \sum_{t=s}^e
\bigl(f_t - \bar{f}_{s,e,t}^d
\bigr)^2 - \sum_{t=s}^e
\bigl(f_t - \bar {f}_{s,e,t}^{\eta_{p_0+r}}
\bigr)^2,
\end{equation}
since obviously $\sum_{t=s}^e (f_t - \bar{f}_{s,e,t}^d)^2 \le\sum_{t=s}^e (f_t - \bar{X}_{s,e,t}^d)^2$.
For any $d$, we have an ANOVA-type decomposition
\[
\sum_{t=s}^e \bigl(f_t -
\bar{f}_{s,e,t}^d\bigr)^2 = \sum
_{t=s}^e \bigl(f_t - \bar{f} - \bigl
\langle f, \psi_{s,e}^d\bigr\rangle\psi_{s,e,t}^d
\bigr)^2 = \sum_{t=s}^e
(f_t - \bar{f})^2 - \bigl\langle f, \psi_{s,e}^d
\bigr\rangle^2.
\]
Therefore, the right-hand side of (\ref{eq:stochdet}) reduces to
\begin{eqnarray*}
\bigl\langle f, \psi_{s,e}^{\eta_{p_0+r}}\bigr\rangle^2 -
\bigl\langle f, \psi _{s,e}^d\bigr\rangle^2 & =
& \bigl(\bigl|\bigl\langle f, \psi_{s,e}^{\eta_{p_0+r}}\bigr\rangle\bigr| - \bigl|\bigl
\langle f, \psi _{s,e}^d\bigr\rangle\bigr|\bigr) \bigl(\bigl|\bigl
\langle f, \psi_{s,e}^{\eta_{p_0+r}}\bigr\rangle\bigr| + \bigl|\bigl\langle f, \psi
_{s,e}^d\bigr\rangle\bigr|\bigr)
\\
& \ge& \bigl(\bigl|\bigl\langle f, \psi_{s,e}^{\eta_{p_0+r}}\bigr\rangle\bigr| -
\bigl|\bigl\langle f, \psi_{s,e}^d\bigr\rangle\bigr|\bigr) \bigl|\bigl
\langle f, \psi_{s,e}^{\eta_{p_0+r}}\bigr\rangle\bigr|.
\end{eqnarray*}
Since by Lemma~\ref{lem:near}, $|\langle f, \psi_{s,e}^{\eta
_{p_0+r}}\rangle|$ is a local
maximum, we can invoke Lemma~2 of \citet{cf12}, by which we obtain
\[
\bigl|\bigl\langle f, \psi_{s,e}^{\eta_{p_0+r}}\bigr\rangle\bigr| - \bigl|\bigl\langle
f, \psi _{s,e}^d\bigr\rangle\bigr| \ge C|d -
\eta_{p_0+r}| T^{-1/2} f'_{p_0+r}.
\]
Combining Lemma~1 of \citet{cf12} with the last assertion of Lemma~\ref
{lem:near}, we
obtain $|\langle f, \psi_{s,e}^{\eta_{p_0+r}}\rangle| \ge C \delta_T
T^{-1/2} f'_{p_0+r}$. This finally
yields
\[
\sum_{t=s}^e \bigl(f_t -
\bar{f}_{s,e,t}^d\bigr)^2 - \sum
_{t=s}^e \bigl(f_t - \bar
{f}_{s,e,t}^{\eta_{p_0+r}}\bigr)^2 \ge C |d -
\eta_{p_0+r}| \delta_T \bigl(f'_{p_0+r}
\bigr)^2 / T.
\]
We decompose the left-hand side of (\ref{eq:stochdet}) as
%
%e15 #&#
\begin{eqnarray}
\label{eq:stochdec}&& 2 \sum_{t=s}^e
\varepsilon_t \bigl(\bar{X}_{s,e,t}^d -
\bar{f}_{s,e,t}^{\eta
_{p_0+r}}\bigr)
\nonumber
\\[-8pt]
\\[-8pt]
\nonumber
&&\qquad = 2 \sum
_{t=s}^e \varepsilon_t \bigl(
\bar{X}_{s,e,t}^d - \bar{f}_{s,e,t}^d
\bigr) + 2 \sum_{t=s}^e
\varepsilon_t \bigl(\bar{f}_{s,e,t}^d -
\bar{f}_{s,e,t}^{\eta
_{p_0+r}}\bigr).
\end{eqnarray}
Without loss of generality, assume $d \ge\eta_{p_0+r}$. The second
term on the right-hand side of (\ref{eq:stochdec}) decomposes as
\begin{eqnarray*}
\sum_{t=s}^e \varepsilon_t
\bigl(\bar{f}_{s,e,t}^d - \bar{f}_{s,e,t}^{\eta
_{p_0+r}}
\bigr) &=& \Biggl( \sum_{t=s}^{\eta_{p_0+r}} + \sum
_{t=\eta_{p_0+r}+1}^d + \sum
_{t=d+1}^e \Biggr) \varepsilon_t \bigl(
\bar{f}_{s,e,t}^d - \bar {f}_{s,e,t}^{\eta_{p_0+r}}
\bigr) \\
&=& I + \mathit{II} + \mathit{III}.
\end{eqnarray*}
We bound
\begin{eqnarray*}
|I| & \le& \sqrt{\eta_{p_0+r}-s+1} \Biggl\llvert \frac{1}{\sqrt{\eta
_{p_0+r}-s+1}}\sum
_{t=s}^{\eta_{p_0+r}}\varepsilon_t\Biggr
\rrvert \\
&&{}\times\Biggl\llvert \frac{1}{d-s+1}\sum_{t=s}^d
f_t - \frac{1}{\eta_{p_0+r}-s+1}\sum_{t=s}^{\eta_{p_0+r}}
f_t\Biggr\rrvert
\\
& \le& \sqrt{\eta_{p_0+r}-s+1} \lambda_2 \frac{C  |d - \eta
_{p_0+r}| f'_{p_0+r}}{\eta_{p_0+r}-s+1}
\le C \lambda_2 |d - \eta_{p_0+r}| f'_{p_0+r}
\delta_T^{-1/2},
\end{eqnarray*}
and we note that the bound for $\mathit{III}$ is of the same order. Similarly,
the bound for $\mathit{II}$ is
$C \lambda_2 |d - \eta_{p_0+r}|^{1/2} f'_{p_0+r}$. The first term on
the right-hand side of (\ref{eq:stochdec}) decomposes as
\[
\sum_{t=s}^e \varepsilon_t
\bigl(\bar{X}_{s,e,t}^d - \bar{f}_{s,e,t}^d
\bigr) = \Biggl(\sum_{t=s}^d + \sum
_{t=d+1}^e \Biggr) \varepsilon_t
\bigl(\bar {X}_{s,e,t}^d - \bar{f}_{s,e,t}^d
\bigr) = \mathit{IV} + V.
\]
Note that $\mathit{IV}$ and $V$ are of the same order. We have
\[
\mathit{IV} = \frac{1}{d-s+1} \Biggl(\sum_{t=s}^d
\varepsilon_t \Biggr)^2 \le \lambda_2^2.
\]

Combining all of the above bounds, there exists a constant $C$ such
that (\ref{eq:stochdet}) holds if
%
%e16 #&#
\begin{eqnarray}
\label{eq:ineq3} &&|d - \eta_{p_0+r}| \delta_T T^{-1}
\bigl(f'_{p_0+r}\bigr)^2
\nonumber
\\[-8pt]
\\[-8pt]
\nonumber
&&\qquad \ge C\max\bigl(
\lambda_2 |d - \eta_{p_0+r}| \delta_T^{-1/2}
f'_{p_0+r}, \lambda_2 |d -
\eta_{p_0+r}|^{1/2} f'_{p_0+r},
\lambda_2^2\bigr).
\end{eqnarray}
These three inequalities yield, respectively, $\delta_T \ge(C \lambda
_2 T / f'_{p_0+r})^{2/3}$,\break $|d - \eta_{p_0+r}| \ge
C^2 \lambda_2^2 T^2 (\delta_T f'_{p_0+r})^{-2}$,  $|d - \eta_{p_0+r}|
\ge C \lambda_2^2 T \delta_T^{-1} (f'_{p_0+r})^{-2}$.\vspace*{1pt}
The last inequality
can be ignored as it is implied by the second if $C \ge1$. The first
inequality can also be ignored as the second
inequality and (\ref{eq:restr}) together imply
\[
C^2 \lambda_2^2 T^2
\delta_T^{-2} \bigl(f'_{p_0+r}
\bigr)^{-2} \le C_1 \lambda_1 T^{1/2}
\bigl(f'_{p_0+r}\bigr)^{-1},
\]
which leads to $\delta_T \ge C C_1^{-1/2} \lambda_2 \lambda_1^{-1/2}
T^{3/4} (f'_{p_0+r})^{-1/2}$,
a stronger requirement that
in the first inequality since $\varpi< 1/2$, but automatically
satisfied since $3/4 + \varpi/2 < \Theta$.
Therefore by the second inequality, $\epsilon_T$ can be taken to be
$\max(1, C^2) \lambda_2^2 T^2 \delta_T^{-2} (f'_{p_0+r})^{-2}$.
It remains for us to note that (\ref{eq:restr}) is automatically
satisfied, as required.
This completes the proof of the lemma.
\end{pf}

%le6.4 #&#
\begin{lemma}
\label{lem:howlarge}
Let $X_t$ follow model (\ref{eq:mod}), and let the assumptions of
Theorem~\ref{th:sbs}
hold. Let $\tilde{X}_{s,e}^b$ be defined by
(\ref{eq:xtilde}).
Assume (\ref{eq:cond0}), (\ref{eq:cond1}), and (\ref{eq:cond2}).
On the event $B_T$ from Lemma~\ref{lem:bonf}, we have
$|\tilde{X}_{s,e}^b| > C T^{\Theta- 1/2 - \varpi}$, where $b = \arg\max_{t: s\le t<e} |\tilde{X}_{s,e}^t|$.
\end{lemma}

\begin{pf} Let $r$ be as in Lemma~\ref{lem:near}. We have
\[
\bigl|\tilde{X}_{s,e}^b\bigr| \ge\bigl|\tilde{X}_{s,e}^{\eta_{p_0+r}}\bigr|
\ge\bigl|\tilde {f}_{s,e}^{\eta_{p_0+r}}\bigr| - \lambda_1 \ge C
\delta_T T^{-1/2}\underline{f}_T -
\lambda_1 > C_1 T^{\Theta- 1/2 - \varpi},
\]
which completes the proof.
\end{pf}

%le6.5 #&#
\begin{lemma}
\label{lem:howsmall}
Let $X_t$ follow model (\ref{eq:mod}), and let the assumptions of
Theorem~\ref{th:sbs}
hold. Let $\tilde{X}_{s,e}^b$ be defined by
(\ref{eq:xtilde}).
For some positive constants $C$, $C'$, let $s, e$ satisfy one of three
conditions:
\begin{longlist}[(iii)]
\item[(i)]
$\exists ! 1 \le p \le N$ such that $s \le\eta_p \le e$ and $(\eta_p
- s + 1) \land(e - \eta_p) \le C\epsilon_T$, or
\item[(ii)]
$\exists  1 \le p \le N$ such that $s \le\eta_p \le\eta_{p+1} \le
e$ and $(\eta_p - s + 1) \lor(e - \eta_{p+1}) \le
C' \epsilon_T$, or
\item[(iii)]
$\exists  1 \le p \le N$ such that $\eta_p < s < e \le\eta_{p+1}$.
\end{longlist}
On the event $B_T \cap A_T$ from Lemma~\ref{lem:vnear}, we have $|\tilde
{X}_{s,e}^b| < C\lambda_2 T^{1-\Theta} + \lambda_1$, where $b = \arg\max_{t: s\le t<e} |\tilde{X}_{s,e}^t|$.
\end{lemma}

\begin{pf} We show case (ii); the remaining two cases are similar and simpler.
\[
\bigl|\tilde{X}_{s,e}^b\bigr| \le\bigl|\tilde{f}_{s,e}^b\bigr|
+ \lambda_1 \le\max\bigl(\bigl|\tilde{f}_{s,e}^{\eta_p}\bigr|,
\bigl|\tilde{f}_{s,e}^{\eta_{p+1}}\bigr|\bigr) + \lambda_1 = \bigl|
\tilde{f}_{s,e}^{\eta_{p_0}}\bigr| + \lambda_1 \le C
\epsilon_T^{1/2} f'_{p_0} +
\lambda_1,
\]
where the last inequality uses the definition of $\tilde{f}_{s,e}^t$.
Continuing, for large $T$,
\[
C\epsilon_T^{1/2}f'_{p_0} +
\lambda_1 \le C\lambda_2 T \delta_T^{-1}
+ \lambda_1 \le C\lambda_2 T^{1-\Theta} +
\lambda_1,
\]
which completes the proof.
\end{pf}

With the use of Lemmas \ref{lem:bonf} to \ref{lem:howsmall}, the proof
of the theorem is simple; the following
occurs on the event $B_T \cap A_T$,
which has probability $\ge1 - C_1 T^{-1}$.
At the start of the algorithm, as $s = 0$ and $e = T-1$,
all conditions for Lemma~\ref{lem:vnear} are met and it finds a
change-point within the distance of
$C\epsilon_T$ from the true change-point, by Lemma~\ref{lem:howlarge}.
Under the assumption of the theorem,
both (\ref{eq:cond1}) and (\ref{eq:cond2}) are satisfied within each
segment until
every change-point in $f_t$ has been identified. Then one of the three
conditions, (i), (ii), or (iii) of Lemma~\ref{lem:howsmall}, are met,
and no further
change-points are detected.
\end{pf*}

\begin{pf*}{Proof of Theorem~\ref{th:wbs}}
We start by defining intervals $\mathcal{I}_i$ between change-points in
such a way that their
lengths are at least of order $\delta_T$, and they are separated from
the change-points
also by distances at least of order $\delta_T$. To fix ideas, define
$\mathcal{I}_i =  [\eta_{i-1}+\frac{1}{3}(\eta_i - \eta_{i-1}), \eta
_{i-1}+\frac{2}{3}(\eta_i - \eta_{i-1}) ]$,
$i = 1, \ldots, N+1$.

Each stage of the algorithm uses CUSUM statistics computed over $M$ intervals
$(s_m, e_m)$, $m = 1, \ldots, M$, drawn uniformly (independently with
replacement) from the set
$\{(s,e) \dvtx s < e, 1 \le s \le T-1, 2 \le e \le T\}$. Define
the event $D_T^M$ as follows:
\[
D_T^M = \bigl\{ \forall i=1,\ldots,N\ \exists m=1,
\ldots,M\ (s_m, e_m) \in \mathcal{I}_i
\times\mathcal{I}_{i+1} \bigr\}.
\]
Note that
\[
P\bigl(\bigl(D_T^M\bigr)^c\bigr) \le\sum
_{i=1}^N \prod
_{m=1}^M \bigl(1 - P\bigl((s_m,
e_m) \in \mathcal{I}_i \times\mathcal{I}_{i+1}
\bigr)\bigr) \le T \delta_T^{-1} \bigl(1 -
\delta_T^2 T^{-2} / 9\bigr)^M.
\]
The remaining arguments will be valid on the set $D_T^M$. If an
interval $(s_m, e_m)$ is such that
$(s_m, e_m) \in\mathcal{I}_i \times\mathcal{I}_{i+1}$, and thus
$(s_m, e_m)$ contains one
change-point only, $\eta_{i}$, then arguing as in \citet{cf12}, Lemma~1,
we have
%
%e17 #&#
\begin{equation}
\label{eq:lem1star} \bigl|\tilde{f}_{s_m,e_m}^{\eta_i}\bigr| = \max
_{t: s_m \le t < e_m} \bigl|\tilde {f}_{s_m, e_m}^t\bigr| \ge C
\delta_T^{1/2} f'_{i}.
\end{equation}

Let $(s, e)$ now be a generic interval satisfying (\ref{eq:cond1}) and
(\ref{eq:cond2}), with $\epsilon_T$
and $\delta_T$ as in the statement of this theorem. The remaining
arguments are valid on the set $B_T$. Consider
%
%e18 #&#
\begin{equation}
\label{eq:amax} (m_0, b) = \arg\max_{(m,t) \dvtx m\in\mathcal{M}_{s,e}, s_m \le t < e_m} \bigl|
\tilde{X}_{s_m, e_m}^t\bigr|,
\end{equation}
where $\mathcal{M}_{s,e} = \{ m \dvtx (s_m, e_m) \subseteq(s, e),  1
\le m \le M \}$.
Imposing the condition
%
%e19 #&#
\begin{equation}
\label{eq:dve} \delta_T \ge3 \epsilon_T,
\end{equation}
we guarantee that both $s$ and $e$ are sufficiently bounded away from
all the previously undetected
change-points $\eta_i \in(s,e)$ in the sense that $\mathcal{I}_i \cup
\mathcal{I}_{i+1} \subset(s,e)$ for
all such $i$. Denote the set of these $i$'s by $\mathcal{J}_{s,e}$. For each
$i \in\mathcal{J}_{s,e}$, there exists an $m_i \in\mathcal{M}_{s,e}$
such that $(s_{m_i}, e_{m_i}) \in\mathcal{I}_i \times\mathcal
{I}_{i+1}$, and thus
%
%e20 #&#
\begin{eqnarray}
\label{eq:howlarge2} \bigl|\tilde{X}_{s_{m_0}, e_{m_0}}^b\bigr| &\ge&\max
_{t: s_{m_i} \le t < e_{m_i}} \bigl|\tilde{X}_{s_{m_i}, e_{m_i}}^t\bigr|
\nonumber
\\[-8pt]
\\[-8pt]
\nonumber
&\ge&\bigl|
\tilde{X}_{s_{m_i}, e_{m_i}}^{\eta_i}\bigr| \ge\bigl|\tilde{f}_{s_{m_i},
e_{m_i}}^{\eta_i}\bigr|
- \lambda_1 \ge C_1 \delta_T^{1/2}
f'_i,
\end{eqnarray}
provided that
%
%e21 #&#
\begin{equation}
\label{eq:dvl1} \delta_T \ge C_8 \lambda_1^2
(\underline{f}_T)^{-2}.
\end{equation}
Therefore,
%
%e22 #&#
\begin{equation}
\label{eq:bfb} \bigl|\tilde{f}_{s_{m_0}, e_{m_0}}^b\bigr| \ge\bigl|
\tilde{X}_{s_{m_0}, e_{m_0}}^b\bigr| - \lambda_1 \ge
C_2 \delta_T^{1/2} \max_{i \in\mathcal{J}_{s,e}}
f'_i.
\end{equation}
By Lemma~2.2 in \citet{venkatraman1993}, there exists a change-point
$\eta_{p_0+r}$, immediately to the left or to the
right of $b$ such that
%
%e23 #&#
\begin{equation}
\label{eq:magn} \bigl|\tilde{f}_{s_{m_0}, e_{m_0}}^{\eta_{p_0+r}}\bigr| > \bigl|
\tilde{f}_{s_{m_0},
e_{m_0}}^b\bigr| \ge C_2 \delta_T^{1/2}
\max_{i \in\mathcal{J}_{s,e}} f'_i.
\end{equation}
Now, the following two situations are impossible:
\begin{longlist}[(1)]
\item[(1)]
$(s_{m_0}, e_{m_0})$ contains one change-point only, $\eta_{p_0+r}$,
and either $\eta_{p_0+r} - s_{m_0}$ or
$e_{m_0} - \eta_{p_0+r}$ is not bounded\vspace*{-1pt} from below by $C_3 \delta_T$;
\item[(2)]
$(s_{m_0}, e_{m_0})$ contains two change-points only, say $\eta
_{p_0+r}$ and $\eta_{p_0+r+1}$, and both
$\eta_{p_0+r} - s_{m_0}$ and $e_{m_0} - \eta_{p_0+r+1}$ are not bounded
from below by $C_3 \delta_T$.
\end{longlist}
Indeed, if either situation were true, then using arguments as in Lemma~\ref{lem:howsmall}, we would
obtain that $\max_{t: s_{m_0} \le t < e_{m_0}} |\tilde{f}_{s_{m_0},
e_{m_0}}^t|$ were not bounded from
below by $C_2 \delta_T^{1/2} \max_{i \in\mathcal{J}_{s,e}} f'_i$,
a contradiction to (\ref{eq:bfb}). This proves that the interval
$(s_{m_0}, e_{m_0})$
satisfies condition (\ref{eq:cond1}) (with $\delta_T$ as in the
statement of this theorem), and thus
we can follow the argument from the proof of Lemma~2 in \citet{cf12} to
establish that if $|b' - \eta_{p_0+r}| = C \gamma_T$ for a certain $C$,
with $\gamma_T = \delta_T^{1/2}\lambda_1 / f'_{p_0+r}$,
and if $\tilde{f}_{s_{m_0}, e_{m_0}}^{\eta_{p_0+r}} > \tilde
{f}_{s_{m_0}, e_{m_0}}^{b'}$
(assuming w.l.o.g. $\tilde{f}_{s_{m_0}, e_{m_0}}^{\eta_{p_0+r}} > 0$), then
$\tilde{f}_{s_{m_0}, e_{m_0}}^{\eta_{p_0+r}} \ge\tilde{f}_{s_{m_0},
e_{m_0}}^{b'} + 2\lambda_1$.

With this result, it is then straightforward to proceed like in the
proof of Lemma~\ref{lem:near} to show that
$|b - \eta_{p_0+r}| \le C_4\gamma_T$, and that $|\tilde
{f}_{s_{m_0},e_{m_0}}^t|$ must have a local
maximum at $t = \eta_{p_0+r}$.

To establish that $|b - \eta_{p_0+r}| \le C_7\epsilon_T$, we need to
use the above results to obtain
an improved version of Lemma~\ref{lem:vnear}.
The arguments in the remainder of the proof are valid on the set $A_T$.
Following the proof of Lemma~\ref{lem:vnear}
for the interval $(s_{m_0}, e_{m_0})$ with $\gamma_T = \delta
_T^{1/2}\lambda_1 / f'_{p_0+r}$, in the notation of
that lemma and using an argument like in Lemma~2 of \citet{cf12}, we obtain
\[
\bigl|\bigl\langle f, \psi_{s_{m_0},e_{m_0}}^{\eta_{p_0+r}}\bigr\rangle\bigr| - \bigl|\bigl\langle
f, \psi_{s_{m_0},e_{m_0}}^d\bigr\rangle\bigr| \ge C|d -
\eta_{p_0+r}| \delta_T^{-1/2} f'_{p_0+r}.
\]
Additionally, by (\ref{eq:magn}), $|\langle f, \psi
_{s_{m_0},e_{m_0}}^{\eta_{p_0+r}}\rangle| \ge
C_5 \delta_T^{1/2} f'_{p_0+r}$, which combined yields
\[
\sum_{t=s_{m_0}}^{e_{m_0}} \bigl(f_t -
\bar{f}_{s_{m_0},e_{m_0},t}^d\bigr)^2 - \sum
_{t=s_{m_0}}^{e_{m_0}} \bigl(f_t -
\bar{f}_{s_{m_0},e_{m_0},t}^{\eta_{p_0+r}}\bigr)^2 \ge C |d - \eta
_{p_0+r}| \bigl(f'_{p_0+r}\bigr)^2.
\]
This in turn leads to the following replacement for the triple
inequality (\ref{eq:ineq3}):
\begin{eqnarray*}
&&|d - \eta_{p_0+r}| \bigl(f'_{p_0+r}
\bigr)^2 \\
&&\qquad\ge C\max\bigl(\lambda_2 |d - \eta
_{p_0+r}| \delta_T^{-1/2}f'_{p_0+r},
\lambda_2 |d - \eta_{p_0+r}|^{1/2}
f'_{p_0+r}, \lambda_2^2\bigr).
\end{eqnarray*}
These three inequalities yield, respectively, $\delta_T \ge C^2 \lambda
_2^2 / (f'_{p_0+r})^2$,
$|d - \eta_{p_0+r}| \ge C^2 \lambda_2^2 / (f'_{p_0+r})^2$, $|d - \eta
_{p_0+r}| \ge C \lambda_2^2 / (f'_{p_0+r})^2$.
The second
inequality and the requirement that $|d - \eta_{p_0+r}| \le C_6 \gamma
_T = C_6 \lambda_1 \delta_T^{1/2} / f'_{p_0+r}$
(see the proof of Lem\-ma~\ref{lem:vnear})
together imply
$\delta_T \ge C^4 C_6^{-2} \lambda_1^{-2} \lambda_2^4
(f'_{p_0+r})^{-2}$. Combining this with the first
inequality, we obtain
%
%e24 #&#
\begin{equation}
\label{eq:lbb} \delta_T \ge C^2 \lambda_2^2
\bigl(f'_{p_0+r}\bigr)^{-2} \max
\bigl(C^2 C_6^{-2} \lambda _1^{-2}
\lambda_2^2, 1\bigr).
\end{equation}
By the second and third inequalities, $\epsilon_T$ can be taken to be
$\max(1, C^2) \lambda_2^2 / \break (f'_{p_0+r})^2$.
At this point, we recall the constraints (\ref{eq:dve}) and (\ref
{eq:dvl1}). Taking $\lambda_1$
and $\lambda_2$ to be of the lowest permissible order of magnitude,
that is $\sim\log^{1/2}T$,
these constraints together with (\ref{eq:lbb}) stipulate that we must
have $\delta_T \ge C_9 \log T / (\underline{f}_T)^2$
for a large enough $C_9$.

With the use of the above results, the proof of the theorem proceeds as
follows; the following
occurs on the event $B_T \cap A_T \cap D_T^M$,
which has probability
$\ge1 - C_1 T^{-1} - T \delta_T^{-1} (1 - \delta_T^2 T^{-2} / 9)^M$.
At the start of the algorithm, as $s = 0$ and $e = T-1$,
(\ref{eq:cond1}) and (\ref{eq:cond2}) (with $\delta_T$ and $\epsilon_T$
as in the statement of this
theorem) are satisfied, and therefore, by formula (\ref{eq:howlarge2}),
the algorithm detects a change-point $b$
on that interval, defined by formula (\ref{eq:amax}). By the above
discussion, $b$ is within the distance
of $C\epsilon_T$ from the change-point. Then (\ref{eq:cond1}) and (\ref
{eq:cond2}) (with $\delta_T$ and
$\epsilon_T$ as in the statement of this theorem) are satisfied within
each segment until
every change-point in $f_t$ has been identified. Once this has
happened, we note that every subsequent
interval $(s_m, e_m)$ satisfies the assumptions on $(s, e)$ from
Lemma~\ref{lem:howsmall} and therefore $|\tilde{X}_{s_{m_0},
e_{m_0}}^b| < C \lambda_2 + \lambda_1 \le\zeta_T$,
which means that no further change-points are detected.
\end{pf*}

\begin{pf*}{Proof of Theorem~\ref{th:ssic}}
The following considerations are valid on the set $A_T \cap B_T \cap
D_T^M$ (from
Theorem~\ref{th:wbs}) which has probability
$\ge1 - C_1 T^{-1} - T \delta_T^{-1} (1 - \delta_T^2 T^{-2} / 9)^M$.
First consider the case $k > N$. Let $\bar{X}_{s,e}$ be the sample mean of
$X_t$ on the interval $[s,e]$ and recall the definition of $\psi
_{s,e}^d$ from
Lemma~\ref{lem:vnear}. The difference $\hat{\sigma}_{k-1}^2 - \hat
{\sigma}_{k}^2$ must
necessarily be of the form
%
%e25 #&#
\begin{eqnarray}\label{eq:ipdv}
\hat{\sigma}^2_{k-1} - \hat{\sigma}^2_{k}
& = & \frac{1}{T} \Biggl\{ \sum_{i=s}^e
(X_i - \bar{X}_{s,e})^2 - \sum
_{i=s}^e \bigl(X_i -
\bar{X}_{s,e} - \bigl\langle X, \psi_{s,e}^d\bigr
\rangle\psi _{s,e,i}^d\bigr)^2 \Biggr\}
\nonumber
\\
& = & \frac{1}{T} \Biggl\{ 2 \sum_{i=s}^e
(X_i - \bar{X}_{s,e}) \bigl\langle X, \psi_{s,e}^d
\bigr\rangle\psi_{s,e,i}^d - \sum
_{i=s}^e \bigl\langle X, \psi_{s,e}^d
\bigr\rangle^2 \bigl(\psi_{s,e,i}^d
\bigr)^2 \Biggr\}
\\
 & = & \frac{\langle X, \psi_{s,e}^d\rangle^2}{T}.\nonumber
\end{eqnarray}
From the proof of Theorem~\ref{th:wbs}, in the case $k > N$, that is,
once all
the change-points have been detected, we have $\langle X, \psi
_{s,e}^d\rangle^2 \le C(\lambda_1^2 + \lambda_2^2)
\le C \log T$. Therefore, for a constant $\upsilon> 0$, and using
the fact that on the set $A_T$, we have $|\hat{\sigma}^2_N - \operatorname
{Var}(\varepsilon_t)| \le C T^{-1}\log T$,
we obtain
\begin{eqnarray*}
\operatorname{sSIC}(k) - \operatorname{sSIC}(N) & = & \frac{T}{2} \log
\frac{\hat{\sigma
}_k^2}{\hat{\sigma}_{N}^2} + (k - N) \log^\alpha T
\\
& = & \frac{T}{2} \log \biggl( 1 - \frac{\hat{\sigma}_{N}^2 - \hat{\sigma
}_{k}^2}{\hat{\sigma}_{N}^2} \biggr) + (k - N)
\log^\alpha T
\\
& \ge& -\frac{T}{2} (1 + \upsilon) \frac{\hat{\sigma}_{N}^2 - \hat
{\sigma}_{k}^2}{\hat{\sigma}_{N}^2} + (k - N)
\log^\alpha T
\\
& \ge& -C_1 \log T + (k - N) \log^\alpha T,
\end{eqnarray*}
which is guaranteed to be positive for $T$ large enough. Conversely, if
$k < N$, then by
formulae (\ref{eq:ipdv}) and (\ref{eq:howlarge2}), we have
$\hat{\sigma}^2_{k} - \hat{\sigma}^2_{k+1} \ge C \delta_T \underline
{f}_T^2 / T$ and hence
\begin{eqnarray*}
\operatorname{sSIC}(k) - \operatorname{sSIC}(N) & = & \frac{T}{2} \log
\frac{\hat{\sigma
}_k^2}{\hat{\sigma}_{N}^2} + (k - N) \log^\alpha T
\\
& = & \frac{T}{2} \log \biggl( 1 + \frac{\hat{\sigma}_{k}^2 - \hat{\sigma
}_{N}^2}{\hat{\sigma}_{N}^2} \biggr) + (k - N)
\log^\alpha T
\\
& \ge& \frac{T}{2} (1 - \upsilon) \frac{\hat{\sigma}_{k}^2 - \hat
{\sigma}_{N}^2}{\hat{\sigma}_{N}^2} - N
\log^\alpha T
\\
& \ge& C \delta_T \underline{f}_T^2 - N
\log^\alpha T,
\end{eqnarray*}
which is again guaranteed to be positive for $T$ large enough. Hence
for $T$ large enough and on
the set $A_T \cap B_T \cap D_T^M$, $\operatorname{sSIC}(k)$ is necessarily
minimised at $N$ and therefore
$\hat{N} = N$, as required.
\end{pf*}

%s7 #&#
\section{Test models and methods used in the simulation study}
\label{sec:simextra}

In the list below, we provide specifications of the test signals $f_t$ and
standard deviations $\sigma$ of the noise $\varepsilon_t$
used in the simulation study of Section~\ref{sec:ss}, as well as
reasons why
these particular signals were used.

\begin{longlist}[(2)]
\item[(1)]
\texttt{blocks}: length
2048,
change-points at $205, 267, 308, 472, 512, 820, 902,\break   1332, 1557, 1598,
1659$, values
between change-points
$0, 14.64, -3.66, 7.32,\break   -7.32,
10.98, -4.39, 3.29, 19.03, 7.68, 15.37, 0$. Standard deviation of the
noise $\sigma= 10$.
Reason for choice: a standard piecewise-constant test signal widely
analysed in the literature.
\item[(2)]
\texttt{fms}: length 497,
change-points at $139, 226, 243, 300, 309, 333$,
values between change-points
$-0.18, 0.08, 1.07, -0.53, 0.16, -0.69, -0.16$. Standard deviation of
the noise $\sigma= 0.3$.
Reason for choice: a test signal proposed in \citet{fms14}.
\item[(3)]
\texttt{mix}: length 560,
change-points at $11, 21, 41, 61, 91, 121, 161, 201, 251,\break   301, 361,
421, 491$,
values between change-points
$7, -7, 6, -6, 5, -5, 4, -4, 3, \break  -3, 2, -2, 1, -1$.
Standard deviation of the noise $\sigma= 4$. Reason for choice: a mix
of prominent change-points
between short intervals of constancy and less prominent change-points
between longer intervals.
\item[(4)]
\texttt{teeth10}: length 140,
change-points at
$11, 21, 31, 41, 51, 61, 71, 81, 91, \break  101, 111, 121, 131$,
values between change-points
$0, 1, 0, 1, 0, 1, 0, 1, 0, 1, 0, 1,  0, 1$.
Standard deviation of the noise $\sigma= 0.4$.
Reason for choice: frequent change-points, occurring every 10th
observation, in the shape of `teeth'.
\item[(5)]
\texttt{stairs10}: length 150,
change-points at
$11, 21, 31, 41, 51, 61, 71, 81, 91,\break   101, 111, 121, 131, 141$,
values between change-points
$1, 2, 3, 4, 5, 6, 7, 8, 9, 10,\break   11, 12, 13, 14, 15$.
Standard deviation of the noise $\sigma= 0.3$.
Reason for choice: frequent change-points, occurring every 10th
observation, in the shape of `stairs'.
\end{longlist}

The list below provides extra details of the competing methods used in the
simulation study.

\begin{longlist}[\quad:]
\item[\texttt{strucchange}:] the main routine for estimating the number
and locations
of change-points is \texttt{breakpoints}. It implements the procedure by
\citet{bp03}.
It is suitable for use in general regression problems, but also in the
signal plus
noise set-up. Given an input vector \texttt{x}, the command we use is
\texttt{breakpoints(x \mbox{$\sim$} 1)}. The \texttt{breakpoints} routine requires a
minimum segment
size, which makes it not fully automatic. The results reported in the
paper are with
the default minimum segment size, which may not be the optimal choice
for our test
signals. We tried changing the minimum segment size to 1, but this
resulted in
execution times that were too long to permit inclusion of the method in
our simulation
study. We refer to the method as `B\&P' throughout the paper.

\item[\texttt{Segmentor3IsBack}:] the main routine is \texttt{Segmentor}. It
implements
a fast algorithm for minimising the least-squares cost function for change-point
detection, as described in \citet{r10}. The function \texttt{SelectModel} then selects the best
model according to (by default) the `oracle' penalisation as described in
\citet{l05}. Our execution is
\begin{verbatim}
z <- Segmentor(x, model=2)
SelectModel(z)
\end{verbatim}
The routine \texttt{Segmentor} requires specification of the maximum
number of
segments, which is set to 15 by default. We do not change this default setting.
None of our test signals exceed this maximum number of segments. We
refer to this
method as `S3IB'.

%f5 #&#
\begin{figure}

\includegraphics{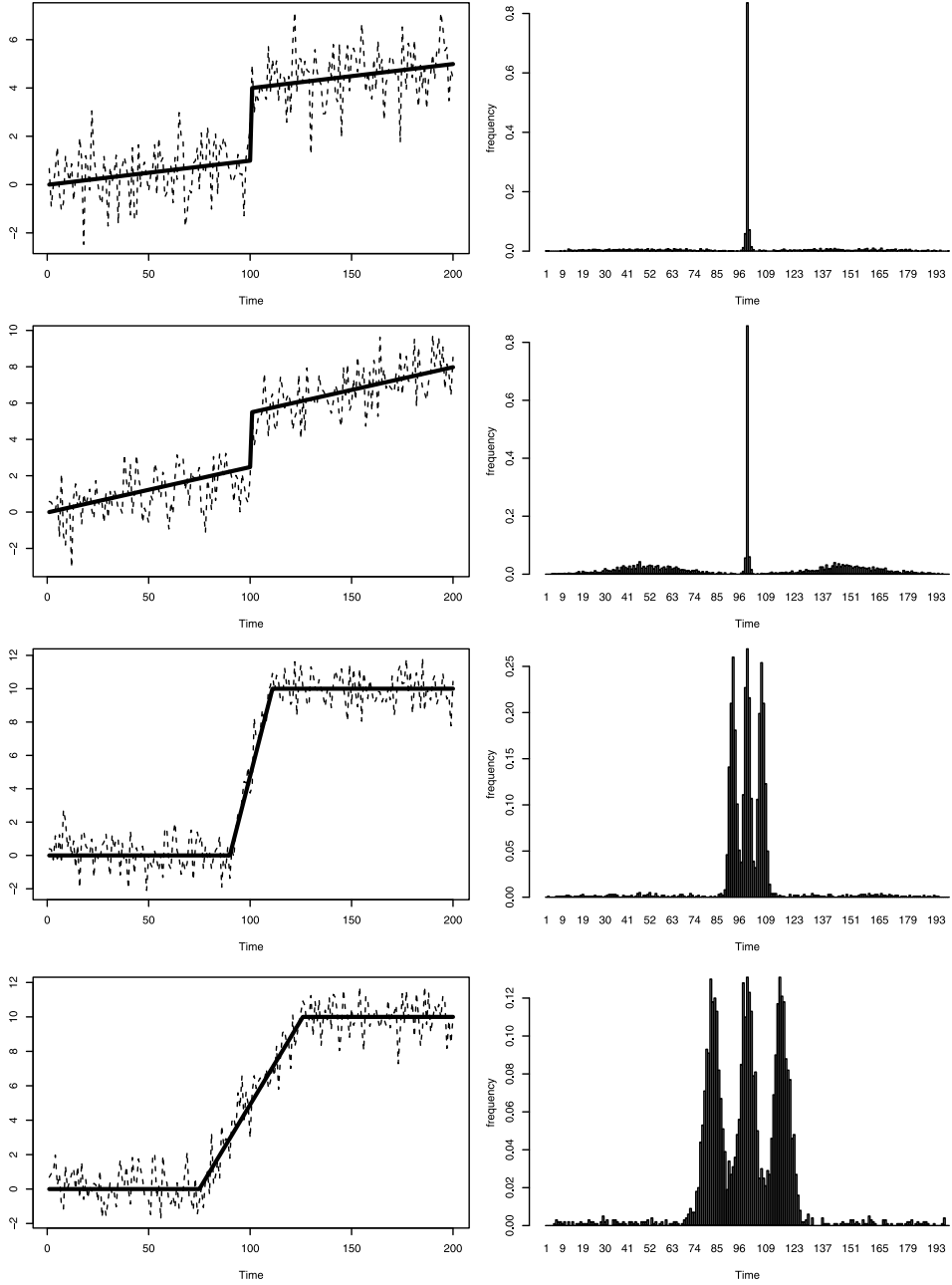}

\caption{Left column: functions $f_t$ containing linear trends
(thick solid lines) and typical realisations of model (\protect\ref
{eq:model}) with
$\varepsilon_t$ i.i.d. standard normal (thin dashed lines).
Right column: the corresponding bar plots of
the frequencies with which change-points were detected at
each time $t$, using the WBS method with threshold constant
$C = 1.3$, over 1000 realisations of each model.}\label{fig:linear}
\end{figure}

\item[\texttt{changepoint}:] the main routine is \texttt{cpt.mean}. It implements
a (different) fast algorithm for minimising the least-squares cost
function for
change-point detection, as described in \citet{kfe12}. The best model is
then selected,
by default, via the SIC penalty. Our
execution is
\begin{verbatim}
cpt.mean(x/mad(diff(x)/sqrt(2)), method="PELT")@cpts,
\end{verbatim}
where the \texttt{mad} function implements the median absolute deviation
estimator.
We refer to this method as `PELT'.

\item[\texttt{cumSeg}:] the main routine is \texttt{jumpoints}, implementing an
algorithm described in \citet{ma11}. We do not change the default
setting which requires
`the starting number of changepoints', which `should be quite larger than
the supposed number of (true) changepoints' (quotes from the package
manual) and is set to
\texttt{min(30, round(length(x)/10))} by default. None of our test signals
violates this. Our execution is \texttt{jumpoints(x)}. We refer to this method
as `cumSeg'.

\item[\texttt{stepR}:] the main routine is \texttt{smuceR}, implementing
a multiscale algorithm described in \citet{fms14}. We leave the default
settings unchanged. Our execution is
\begin{verbatim}
smuceR(x, 1:length(x), family="gauss")
\end{verbatim}
We refer to this method as `SMUCE'.
\end{longlist}

%s8 #&#
\section{Performance of WBS in the presence of linear trends}
\label{sec:linear}

Figure~\ref{fig:linear} shows the results of a small-scale
simulation study aimed at obtaining some insight into the
performance of WBS under model misspecification, namely in
cases where the true function $f_t$ exhibits linear
trends.

In the example from the top row of that figure, the linear trends are so
flat that they are almost completely ignored by WBS. However,
in the example from the second row, the linear trends are more
pronounced, and spurious detection of change-points within
the trend sections tends to occur towards their middle
parts. This can be interpreted in at least two ways: (i) WBS
considers the middle part of a section with a linear trend as the most
likely location of a change-point in the piecewise-constant approximation
of that linear trend, which is natural, and (ii) the change-points
(spuriously) detected
within the trend sections tend to be separated from the main (correctly
detected)
change-point in the middle of the time domain, which is beneficial
for the interpretability of the main change-point.

In the bottom two examples,
spurious detection of change-points within the trend sections tends to
occur towards their middle parts and towards their edges. This can be
interpreted
as the algorithm producing piecewise-constant approximations to the linear
trends in which the change-points are spaced out rather than being
clustered together, which hopefully leads to those approximations being visually
attractive.
\end{appendix}

% imsref loaded by akundreckaite, 2014-07-23 09:47:16
%


%suskaldyti doi

\printaddresses
\end{document}